%% file: ContLu2014.tex
\documentclass[10pt]{article} % use larger type; default would be 10pt
 
\usepackage[utf8]{inputenc} % set input encoding (not needed with XeLaTeX)

\usepackage{geometry} % to change the page dimensions
\usepackage{graphicx} % support the \includegraphics command and options

 \usepackage{array} % for better arrays (eg matrices) in maths
\usepackage{mathtools}
\usepackage{amsfonts,hyperref,url}
\usepackage{amssymb}
\usepackage{amsthm}

\newtheorem{hyp}{Assumption}
\newtheorem{lemma}{Lemma}[section]
\newtheorem{theorem}{Theorem}[section]
\newtheorem{rmq}{Remark}[section]
\newtheorem{example}{Example}[section]
\newtheorem{cor}{Corollary}[section]
\newtheorem{defi}{Definition}
\newtheorem{prop}{Proposition}[section]

\newcommand{\conc}[1]{\underset{#1}{\oplus}}

\everymath{\displaystyle}

\begin{document}
 
\title{Weak approximation of martingale representations}
\author{Rama CONT and Yi LU}

\date{November 2014. Revision: September 2015.} % Activate to display a given date or no date (if empty),
         % otherwise the current date is printed

\maketitle
\begin{abstract}
We present a systematic method for computing explicit approximations to martingale representations for a large class of Brownian functionals.
The approximations are obtained by obtained by computing a directional derivative of the weak Euler scheme and yield a  consistent estimator for the integrand in the martingale representation formula for any square-integrable functional of the solution of an SDE with path-dependent coefficients.  Explicit convergence rates  are derived for functionals which are Lipschitz-continuous in the supremum norm.
Our results require neither the Markov property, nor any differentiability conditions on the functional or the coefficients of the stochastic differential equations involved.
\end{abstract}
\tableofcontents
\newpage

\input WeakApproximation.tex

\input SPA2851.bbl

\end{document}

%% file: WeakApproximation.tex
\section{Introduction}

Let $ W $ be a standard $ d $-dimensional Brownian motion defined on a probability space $ ( \Omega, \mathcal{F}, \mathbb{P} ) $ and $ (\mathcal{F}_t) $ its ($ \mathbb{P} $-completed) natural filtration.
  Then, for any square-integrable $ \mathcal{F}_T $-measurable random variable $ H $, or equivalently, any square-integrable $ (\mathcal{F}_t) $-martingale $ Y (t) = \mathbb{E} [ H | \mathcal{F}_t ] $, there exists a unique $ (\mathcal{F}_t) $-predictable process $ \phi $ with $ \mathbb{E} \left[ \int_0^T \mathrm{tr} ( \phi(u) {^t \phi(u)}) du  \right] < \infty $ such that:
\begin{equation} H= Y (T) = \mathbb{E} [H] + \int_0^T \phi\cdot d W. \label{eq:martingalerepresentation}\end{equation}
The classical proof of this representation result (see e.g. \cite{RY99}) is non-constructive. However in many applications, such as stochastic control or mathematical finance, one is interested in an explicit expression for $ \phi $, which represents an optimal control or a hedging strategy.
Expressions for the integrand $ \phi $ have been derived using a variety of methods and assumptions, using Markovian   techniques \cite{D80, EK88, FR09, JMP00, PP92}, integration by parts \cite{chenglasserman07} or  Malliavin calculus \cite{B81, C70, H79, KOL91, N09, O84, FLLLT99}. Some of these methods are limited to the case where $ Y $ is a Markov process; others require  differentiability of $H$ in the Fr\'echet or Malliavin sense \cite{C70,O84,FLLLT99,FLLL01} or an explicit form for the density \cite{chenglasserman07}.
Almost all of these methods invariably involve an approximation step, either through the solution of an auxiliary partial differential equation (PDE) or the simulation of an auxiliary stochastic differential equation.

A systematic approach to obtaining martingale representation formulae has been proposed in \cite{CF13}, using the Functional It\^o calculus \cite{D09,CF10A,CF10B}: it was shown in \cite[Theorem 5.9]{CF13} that for any square-integrable $ (\mathcal{F}_t) $-martingale $ Y $,
\[ \forall t \in [0, T], \quad Y(t) = Y(0) + \int_0^t \nabla_W Y \cdot dW \quad \mathbb{P} \mathrm{-a.s.} \]
where $ \nabla_W Y $ is the weak vertical derivative of $ Y $ with respect to $ W $, constructed as an $L^2$ limit of pathwise directional derivatives.
This approach does not rely on any Markov property nor on the Gaussian structure of the Wiener space and is  applicable to functionals of a large class of It\^o processes.

In the present work we build on this approach to propose a general framework for computing explicit approximations to the integrand $\phi$ in a general setting in which $H$ is allowed to be a  functional of the solution of a stochastic differential equation (SDE) with path-dependent coefficients:
\begin{equation} \label{sde}
 dX(t) = b( t, X_t ) dt+ \sigma ( t, X_t ) dW(t) \quad X(0) = x_0 \in \mathbb{R}^d
\end{equation}
where $X_t=X(t\wedge .)$ designates the path stopped at $t$ and $$b: [0,T]\times D( [0, T], \mathbb{R}^d ) \rightarrow  \mathbb{R}^d , \qquad\sigma: [0,T]\times D( [0, T], \mathbb{R}^d ) \rightarrow \mathbb{M}_d( \mathbb{R} ) $$ are continuous non-anticipative  functionals.
For any square-integrable variable of the form  $ H = g ( X(t), 0 \leq t \leq T ) $  where $ g: (D( [0, T], \mathbb{R}^d),\|.\|_\infty) \to \mathbb{R}$ is a continuous functional, we construct an explicit sequence of approximations $ \phi_n $ for the integrand $\phi$ in \eqref{eq:martingalerepresentation}.
These approximations are constructed as vertical derivatives, in the sense of the functional It\^o calculus\cite{CF13}, of the weak Euler approximation $F_n$ of the martingale $Y$, obtained by replacing $X$ by the corresponding Euler scheme ${_n X}$:
$$ \phi_n (t) =  \nabla_{\omega} F_n (t, W),\quad{\rm where}\quad F_n (t, \omega) = \mathbb{E} \left[ g ( {_n X} ( \omega \conc{t} W ) ) \right]$$
 where $\conc{t}$ is the concatenation of paths at $t$ and $\nabla_{\omega} F_n$ is the Dupire derivative \cite{D09,Cont2012}, a directional derivative defined as a pathwise limit of  finite-diference approximations. and thus readily computable path-by-path in a simulation setting. 

The main results of the paper are the following.
We first show the existence and continuity of these pathwise derivatives in Theorem \ref{th_diff}.
The convergence of the approximations  $\phi_n$ to the integrand $\phi$  in 
 \eqref{eq:martingalerepresentation}  is shown in Proposition \ref{cor_as2}. Under a Lipschitz assumption on $g$, we provide in Theorem \label{th_error} an $L^{p}$ error estimate for the approximation error.
The proposed approximations are easy to compute and  readily integrated in commonly used numerical schemes for SDEs.

%Unlike methods developed for Markov processes \cite{D80, EK88, FR09, JMP00, PP92},
Our approach requires neither the Markov property of the underlying processes  nor the differentiability of coefficients, and is thus applicable to functionals of a
large class of semimartingales. By contrast to methods based on Malliavin calculus \cite{B81, C70, H79, KOL91, O84, FLLLT99}, it does not require Malliavin differentiability of the terminal variable $H$ nor does it involve any choice of 'Malliavin weights', a delicate step in these methods.

Ideas based on Functional It\^o calculus have also been recently used by Le{\~a}o  and Ohashi \cite{MR3098445} for weak approximation of Wiener functionals, using a space-filtration discretization scheme. However, unlike the approach proposed in \cite{MR3098445}, our approach is based on a Euler approximation on a fixed time grid, rather than the random time grid used in \cite{MR3098445}, which involves a sequence of first passage times. Our approach  is thus much easier to implement and analyze and is readily integrated in commonly used numerical schemes for approximations of SDEs, which are typically based on fixed time grids.

\paragraph{Outline} We first recall some key concepts and results from the Functional It\^o calculus in section \ref{sec:FunctionalIto}. Section \ref{sec:Euler}  provides some estimates for the path-dependent SDE \eqref{sde} and studies some properties of the Euler approximation for this SDE.
In  Section \ref{sec:approximation} we show that the weak Euler
  approximation (Definition \ref{def:weakEuler}) may be used to approximate
 any square-integrable martingale adapted to the filtration of $X$  by a sequence of smooth functionals of $X$, in the sense of the functional It\^o calculus. Moreover, we provide explicit expressions for the functional derivatives of these approximations.
Section \ref{sec:convergence} analyzes the convergence of this approximation and provides error estimates in Theorem \ref{th_error}. Finally, in Section \ref{sec:comparison} we compare our approximation method with those based on Malliavin calculus.

\noindent \textbf{Notations}: In the sequel, we shall denote by $ \mathcal{M}_{d, n} (\mathbb{R}) $ the set of all $ d \times n $ matrices with real coefficients. We simply denote $ \mathbb{R}^d = \mathcal{M}_{d, 1} (\mathbb{R}) $ and $ \mathcal{M}_d (\mathbb{R}) =   \mathcal{M}_{d, d} (\mathbb{R}) $. For $ A \in \mathcal{M}_d (\mathbb{R}) $, we shall denote by $ ^t A $ the transpose of $ A $, and $ \| A \| = \sqrt{\mathrm{tr} \left( ^t A A \right)} $ the Frobenius norm of $ A $. For $ x, y \in \mathbb{R}^d $, $ x \cdot y $ is the scalar product on $ \mathbb{R}^d $.

\noindent Let $ T > 0 $. We denote by $ D ( [0,T], \mathbb{R}^d ) $ the space of functions defined on $ [0, T] $ with values in $ \mathbb{R}^d $ which are right continuous with left limits (c\`adl\`ag). For a path $ \omega  \in D ( [0,T], \mathbb{R}^d ) $ and $ t \in [0, T] $, we denote by:
\begin{itemize}
\item
$ \omega (t) $ the value of $ \omega $ at time $ t $,
\item
$ \omega (t-) = \lim_{s \to t, s<t} \omega(s) $ its left limit at $ t $,
\item
$ \omega_t = \omega(t \wedge \cdot) $ the path of $ \omega $ stopped at $ t $
\item
$ \omega_{t-} = \omega {\mathbf 1}_{[0, t)} + \omega (t-) {\mathbf 1}_{[t, T]} $
\item
$ \| \omega \|_{\infty} = \sup \{ | \omega(t) |, t \in [0, T] \} $ the supremum norm.
\end{itemize}
We note that $ \omega_t $ and $ \omega_{t-} $ are elements of $ D ( [0,T], \mathbb{R}^d ) $. For a c\`adl\`ag stochastic process $ X $, we shall similarly denote $X_t(.)=X(t\wedge .)$ and $X_{t-}=X {\mathbf 1}_{[0, t)} + X(t-) {\mathbf 1}_{[t, T]}$.

\section{Functional It\^o calculus}\label{sec:FunctionalIto}

\input FIC.tex

\section{Euler approximations for path-dependent SDEs}\label{sec:Euler}

Let $ W $ be a standard d-dimensional Brownian motion defined on a probability space $(\Omega, \mathcal{F}, \mathbb{P} ) $ and $ (\mathcal{F}_t) $ its ($ \mathbb{P} $-completed) natural filtration. We consider the following stochastic differential equation with path-dependent coefficients
\eqref{sde}:
$$
dX(t) = b( t, X_t ) dt+ \sigma ( t, X_t ) dW(t), \qquad X(0) = x_0 \in \mathbb{R}^d $$
where $ b: \Lambda_T \rightarrow \mathbb{R}^d,\
\sigma: \Lambda_T \rightarrow \mathbb{M}_d( \mathbb{R}) $ are non-anticipative  maps, assumed to be Lipschitz continuous with respect to the following distance $ d $ defined on $ \Lambda_T $:
\begin{eqnarray*}
d ((t, \omega), (t', \omega')) = \sup_{u \in [0, T]} | \omega (u \wedge t) - \omega' (u \wedge t') | + \sqrt{| t - t' |}
= \| \omega_t - \omega'_{t'} \|_{\infty} + \sqrt{| t - t' |}
\end{eqnarray*}

\begin{hyp}[Lipschitz continuity of coefficients] \label{hyp_sigma}
$ b: (\Lambda_T,d) \rightarrow \mathbb{R}^d, \sigma: (\Lambda_T, d) \rightarrow \mathbb{M}_d( \mathbb{R} ) $ are Lipschitz continuous: $$ \exists K_{Lip} > 0, \quad \forall t, t' \in [0, T],  \forall \omega, \omega' \in D( [0, T], \mathbb{R}^d ), $$ 
\[ | b(t, \omega) - b(t', \omega') | + \left\| \sigma(t, \omega) - \sigma(t', \omega') \right\| \leq K_{Lip} \: d \left( (t, \omega), (t', \omega') \right). \] 
\end{hyp}

\begin{rmq}
This Lipschitz condition with respect to the distance $ d $ is weaker than a Lipschitz condition with respect to the distance $ d_{\infty} $ introduced in the previous section: it allows for a H\"older smoothness of degree 1/2 in the $t$ variable. 
\end{rmq}
Under  Assumption \ref{hyp_sigma},  \eqref{sde} has a unique strong, $ (\mathcal{F}_t) $-adapted, solution $X$.
\begin{prop} \label{prop_sde}
Under Assumption \ref{hyp_sigma}, there exists a unique $ (\mathcal{F}_t) $-adapted  process $ X $ satisfying \eqref{sde}. Moreover for $ p \geq 1 $, we have:
\begin{equation} \mathbb{E} \left[ \| X_T \|_{\infty}^{2p} \right] \leq C ( 1 + |x_0|^{2p} ) e^{CT} \label{eq:estimate}\end{equation}
for some constant $ C = C (p, T, K_{Lip}) $ depending on $ p , T $ and $ K_{Lip} $.
\end{prop}
\begin{rmq}
Assumption \ref{hyp_sigma} might seem to be quite strong. Indeed, the previous proposition still holds under weaker conditions. For example, the H\"older condition with respect to differences in $t$ can be replaced by the weaker condition: $ \sup_{t \in [0, T]} | b (t, \bar{0} | + \left\| \sigma(t, \bar{0}) \right\| < \infty $ where $ \bar{0} $ denotes the path which takes constant value $ 0 $. However, this assumption is necessary for the convergence of the Euler approximation described later in this section, and especially the results concerning its rate of convergence.
\end{rmq}
\begin{proof}
Existence and uniqueness of a strong solution follows from \cite{P05} (Theorem 7 , Chapter 5):  see \cite[Section 5]{Cont2012}. Let us  prove \eqref{eq:estimate}.
Using the Burkholder-Davis-Gundy inequality and H\"older's inequality, we have:
\begin{eqnarray*}
\mathbb{E} \left[ \| X_T \|_{\infty}^{2p} \right]
& \leq & C(p) \left( \left| x_0 \right|^{2p} + \mathbb{E} \left[ \left( \int_0^T | b (t, X_t) |^2 dt \right)^p \right] + \mathbb{E} \left[ \left( \int_0^T \left\| \sigma (t, X_t) \right\|^2 dt \right)^p \right] \right) \\
& \leq & C(p, T) \left( \left| x_0 \right|^{2p} + \mathbb{E} \left[ \int_0^T | b (t, X_t) |^{2p} dt \right] + \mathbb{E} \left[ \int_0^T \left\| \sigma (t, X_t) \right\|^{2p} dt \right] \right) \\
& \leq & C(p, T) \left( \left| x_0 \right|^{2p} + \mathbb{E} \left[ \int_0^T \left( | b (0, \bar{0}) | + \left\| \sigma (0, \bar{0}) \right\| + K_{Lip} ( \sqrt{t} + \| X_t \|_{\infty} ) \right)^{2p} dt \right] \right) \\
& \leq & C(p, T, K_{Lip}) \left( \left| x_0 \right|^{2p} + 1 + \int_0^T \mathbb{E} \| X_t \|_{\infty}^{2p} dt \right)
\end{eqnarray*}
And we conclude by Gronwall's inequality.
\qedhere
\end{proof}
In the following, we always assume that Assumption \ref{hyp_sigma} holds. The strong solution $ X $ of equation \eqref{sde} is then a semimartingale and defines a non-anticipative functional $X:{\cal W}_T\to \mathbb{R}^d$ given by the It\^o map associated to \eqref{sde}.

\subsection{Euler approximations as  non-anticipative functionals}
We now consider an Euler approximation for the SDE \eqref{sde} and study its properties  as a non-anticipative functional.
Let $ n \in \mathbb{N}, \delta = \frac{T}{n} $. The Euler approximation $ {_n X} $ of $X$  on the grid  $ (t_j = j\delta, j=0..n)$ is defined as follows:
\begin{defi} \label{def_euler}[Euler scheme] For $\omega \in D( [0, T], \mathbb{R}^d ) $, we denote   by  $ {_n X} (\omega)  \in D( [0, T], \mathbb{R}^d ) $ the piecewise constant Euler approximation for \eqref{sde} computed along the path $\omega$, defined   as follows: $ {_n X} (\omega) $ is constant in each interval $ [ t_j, t_{j+1} ) $, $ \forall \: 0 \leq j \leq n-1 $ with $ {_n X} (0, \omega) = x_0 $ and 
\begin{equation} \label{euler}
{_n X} (t_{j+1}, \omega) = {_n X} (t_j, \omega) + b (t_j, {_n X}_{t_j} (\omega) ) \delta + \sigma (t_j, {_n X}_{t_j} (\omega) ) ( \omega (t_{j+1}-) - \omega (t_j-) ), 
\end{equation}
where $ {_n X}_t (\omega)= {_n X} (t\wedge ., \omega) $ and by convention $ \omega (0-) = \omega (0) $.
\end{defi}

When computed along the path of the  Brownian motion $W,$ $ {_n X}(W) $ is simply the  piecewise constant Euler-Maruyama scheme \cite{pardouxtalay1985} for the  stochastic differential equation \eqref{sde}.

By definition, the path $ {_n X} (\omega) $ depends only on a finite number of increments of $ \omega $: $ \omega (t_1-) - \omega (0) $, $ \cdots $, $ \omega (t_n-) - \omega (t_{n-1}-) $. We can thus define a map
$$ p_n: \mathcal{M}_{d, n} (\mathbb{R}) \rightarrow D ([0, T], \mathbb{R}^d) $$ such that for $\omega\in D([0, T]), \mathbb{R}^d)$
\begin{equation} \label{p}
  p_n (\omega (t_1-) - \omega (0) , \omega (t_2-) - \omega (t_{1}-),\cdots, \omega (t_n-) - \omega (t_{n-1}-)) = {_n X} (\omega).
\end{equation}
By a slight abuse of notation,  we denote $ p_t (y) $ the  path  $p_n(y) $ stopped at $ t $.

The map $ p_n: \mathcal{M}_{d, n} (\mathbb{R}) \rightarrow \left( D ([0, T], \mathbb{R}^d), \| \cdot \|_{\infty} \right) $ is then locally Lipschitz continuous, as shown by the following lemma.
\begin{lemma} \label{lem_euler} For every $\eta >0$, there exists a constant $C(\eta,K_{Lip}, T)$ such that for any
 $ y = (y_1, \cdots, y_n),  y' = (y'_1, \cdots, y'_n) \in \mathcal{M}_{d, n} (\mathbb{R}), $
\[  \max_{1 \leq k \leq n} | y_k|\vee | y'_k | \leq \eta \quad \Rightarrow\quad
\| p_n (y) - p_n (y') \|_{\infty} \leq C (\eta, K_{Lip}, T) \max_{1 \leq k \leq n} | y_k - y'_k |. \]
\end{lemma}
\begin{proof}
As the two paths $ p_n (y) $ and $ p_n (y') $ are stepwise constant by construction, it suffices to prove the inequality at times $ (t_j)_{0 \leq j \leq n} $. We prove by induction that:
\begin{equation} \label{rec}
\| p_{t_j} (y) - p_{t_j} (y') \|_{\infty} \leq  C ( \eta, K_{Lip}, T) \max_{1 \leq k \leq j} | y_k - y'_k |
\end{equation}
with some constant $ C $ which depends only on  $ \eta $, $ K_{Lip}, T $ (and $ n $).

For $ j = 0 $, this is clearly the case as $ p (y) (0) = p (y') (0) = x_0 $. Assume that \eqref{rec} is verified for some $ 0 \leq j \leq n-1 $, consider now $ \| p_{t_{j+1}} (y) - p_{t_{j+1}} (y') \|_{\infty} $, we have:
\[ p_n (y) (t_{j+1}) = p_n (y) (t_j) + b (t_j, p_{t_j} (y)) \delta + \sigma (t_j, p_{t_j} (y)) y_{j+1} \]
and
\[ p_n (y') (t_{j+1}) = p_n (y') (t_j) + b (t_j, p_{t_j} (y')) \delta + \sigma (t_j, p_{t_j} (y')) y'_{j+1}. \]
Thus
\begin{eqnarray*}
&& | p_n (y) (t_{j+1}) - p_n (y') (t_{j+1}) | \\
& \leq & | p_n (y) (t_j) - p_n (y') (t_j) | + | b (t_j, p_{t_j} (y)) - b (t_j, p_{t_j} (y')) | \delta \\
&& + \| \sigma (t_j, p_{t_j} (y)) \| \cdot | y_{j+1} - y'_{j+1} | + \| \sigma (t_j, p_{t_j} (y)) - \sigma (t_j, p_{t_j} (y')) \| \cdot | y'_{j+1} | \\
& \leq & C (\eta, K_{Lip}, T) \max_{1 \leq k \leq j} | y_k - y'_k | + K_{Lip} C (\eta, K_{Lip}, T) \max_{1 \leq k \leq j} | y_k - y'_k | \delta \\
&& + \left( \| \sigma (0, \bar{0}) \| + K_{Lip} (\sqrt{t_j} + \| p_{t_j} (y)) \|_{\infty}) \right) | y_{j+1} - y'_{j+1} | + K_{Lip} C (\eta, K_{Lip}, T) \eta \max_{1 \leq k \leq j} | y_k - y'_k |  \\
& \leq & C ( \eta, K_{Lip}, T) \max_{1 \leq k \leq j+1} | y_k - y'_k |
\end{eqnarray*}
(The constant $ C $ may differ from one line to another). \\
And consequently we have:
\[ \| p_{t_{j+1}} (y) - p_{t_{j+1}} (y') \|_{\infty} \leq C ( \eta, K_{Lip}, T) \max_{1 \leq k \leq j+1} | y_k - y'_k | \]
for some different constant $ C $ depending only on  $ \eta $, $ K_{Lip} $ and $ T $ (and $n$). And we conclude by induction.
\qedhere
\end{proof}

\subsection{Strong convergence}
To simplify the notations, $ {_n X}_T (W_T) $ will be noted simply $ {_n X}_T $ in the following.
The following result, which gives a uniform estimate of the discretization error,
 $ X_T - {_n X}_T $ extends similar results known in the Markovian case \cite{F92,pardouxtalay1985,higham02} to the path-dependent SDE \eqref{sde} (see also \cite{MR2409016}):
\begin{prop} \label{prop_euler}
Under \textbf{Assumption \ref{hyp_sigma}} we have the following estimate in  $ L^{2p} $ for the strong error of the piecewise constant Euler-Maruyama scheme:
\[ E\left( \sup_{s \in [0, T]} \| X(s) - {_n X} (s) \|^{2p} \right) \leq C (x_0, p, T, K_{Lip}) \left(\frac{1 + \log n}{n}\right)^p, \quad \forall p \geq 1 \]
with $ C $ a constant depending only on $ x_0 $, $ p $, $ T $ and $ K_{Lip} $.
\end{prop}
\begin{proof}
%The proof of this proposition follows from the same methods used in the Markovian case, see e.g.
The idea is to construct a 'Brownian interpolation' $ {_n \hat{X}}_T $ of  $ {_n X}_T $:
\[ {_n \hat{X}} (s) = x_0 + \int_0^s b \left( \underline{u}, {_n X}_{\underline{u}} \right) du + \int_0^s \sigma \left( \underline{u}, {_n X}_{\underline{u}} \right) dW (u) \]
where $ \underline{u} = \left\lfloor \frac{u}{\delta} \right\rfloor \cdot \delta $ is the largest subdivision point which is smaller or equal to $ u $.

Clearly  $ {_n \hat{X}} $  is a continuous semimartingale and $ \| \sup_{s \in [0, T]} | X (s) - {_n X} (s) | \|_{2p} $ can be controlled by the sum of the two following terms:
\begin{equation} \label{sum}
\| \sup_{s \in [0, T]} | X (s) - {_n X} (s) | \|_{2p} \leq \| \sup_{s \in [0, T]} | X (s) - {_n \hat{X}} (s) | \|_{2p} + \| \sup_{s \in [0, T]} | {_n \hat{X}} (s) - {_n X} (s)| \|_{2p}
\end{equation}
We start with the term $ \| \sup_{s \in [0, T]} | X (s) - {_n \hat{X}} (s) | \|_{2p} $. Using the Burkholder-Davis-Gundy inequality and H\"older's inequality, we have
\begin{eqnarray*}
\mathbb{E} \| X_T - {_n \hat{X}}_T \|_{\infty}^{2p}
& \leq & C(p) \left( \mathbb{E} \left[ \int_0^T | b (s, X_s) - b (\underline{s}, {_n X}_{\underline{s}}) | ds \right]^{2p} + \mathbb{E} \left[ \int_0^T \left\| \sigma (s, X_s) - \sigma (\underline{s}, {_n X}_{\underline{s}}) \right\|^2 ds \right]^p \right) \\
& \leq & C (p, T) \left( \mathbb{E} \left[ \int_0^T | b (s, X_s) - b (\underline{s}, {_n X}_{\underline{s}}) |^{2p} ds \right] + \mathbb{E} \left[ \int_0^T \left\| \sigma (s, X_s) - \sigma (\underline{s}, {_n X}_{\underline{s}}) \right\|^{2p} ds \right] \right) \\
& \leq & C (p, T, K_{Lip}) \: \mathbb{E} \left[ \int_0^T \left( (s - \underline{s})^{p} +  \| X_s - {_n X}_s \|_{\infty}^{2p} \right) ds \right] \\
& \leq & C (p, T, K_{Lip}) \: \left( \frac{1}{n^{p}} +  \int_0^T \mathbb{E} \| X_s - {_n X}_s \|_{\infty}^{2p} \: ds \right)
\end{eqnarray*}
We have used $ {_n X}_{\underline{s}} = {_n X}_s $ as $ {_n X} $ is piecewise constant.

Consider now the second term $ \| \sup_{s \in [0, T]} | {_n \hat{X}} (s) - {_n X} (s) | \|_{2p} $. Noting that:
\[ {_n \hat{X}} (s) - {_n X} (s) = {_n \hat{X}} (s) - {_n \hat{X}} ( \underline{s} ) = b \left( \underline{s}, {_n X}_{\underline{s}} \right) ( s - \underline{s} ) + \sigma \left( \underline{s}, {_n X}_{\underline{s}} \right) ( W (s) - W (\underline{s}) ), \]
we have
\[ \| {_n \hat{X}}_T - {_n X}_T \|_{\infty} \leq C (K_{Lip}, T) ( 1 + \| {_n X}_T \|_{\infty} ) \left( \frac{1}{n} + \sup_{s \in [0, T]} | W (s) - W (\underline{s}) | \right) \]
and
\begin{eqnarray*} 
\mathbb{E} \| {_n \hat{X}}_T - {_n X}_T \|_{\infty}^{2p} 
& \leq & C (p,  K_{Lip}, T) \frac{1}{n^{2p}} \: \mathbb{E} \left[ ( 1 + \| {_n X}_T \|_{\infty} )^{2p} \right] \\
&& +  C (p,  K_{Lip}, T) \: \mathbb{E} \left[ ( 1 + \| {_n X}_T \|_{\infty} ) \sup_{s \in [0, T]} | W (s) - W (\underline{s}) | \right]^{2p}. 
\end{eqnarray*}
By the Cauchy-Schwarz inequality, we have:
\begin{equation} \label{cs}
\mathbb{E} \| {_n \hat{X}}_T - {_n X}_T \|_{\infty}^{2p} \leq C (p, K_{Lip}, T) \left( 1 + \sqrt{\mathbb{E} \| {_n X}_T \|_{\infty}^{4p}} \right) \left( \frac{1}{n^{2p}} + \sqrt{\mathbb{E} \sup_{s \in [0, T]} | W (s) - W (\underline{s}) |^{4p}} \right) 
\end{equation}
We will make use of the following  result \cite[p. 203]{pages2014}:
\[ \forall p > 0, \quad \| \sup_{s \in [0, T]} | W (s) - W (\underline{s}) | \|_p \leq C (W, p) \sqrt{\frac{T}{n} (1+\log n)} \]
which results from the following lemma \cite[Lemma 7.1]{pages2014}:
\begin{lemma} \label{lem_iid}
Let $ Y_1 $, $ \cdots $, $ Y_n $ be non-negative random variables with the same distribution satisfying $ \mathbb{E} \left( e^{\lambda Y_1} \right) < \infty $ for some $ \lambda > 0 $. Then we have:
\[ \forall p > 0, \quad \left\| \max (Y_1, \cdots, Y_n) \right\|_p \leq \frac{1}{\lambda} (\log n + C (p, Y_1, \lambda)) \]
\end{lemma}

We have thus:
\begin{equation} \label{bm}
\sqrt{\mathbb{E} \sup_{s \in [0, T]} | W (s) - W (\underline{s}) |^{4p}} \leq C (p, T) \left( \frac{1 + \log n}{n} \right)^p
\end{equation}
Furthermore, using again the Burkholder-Davis-Gundy inequality, we have:
\begin{eqnarray*}
 \mathbb{E} \| {_n X}_T \|_{\infty}^{4p} & \leq & \mathbb{E} \| {_n \hat{X}}_T \|_{\infty}^{4p} \\
& \leq & C (p) \left( x_0^{4p} + \mathbb{E} \left( \int_0^T | b (\underline{s}, {_n X}_{\underline{s}}) | ds \right)^{4p} + \mathbb{E} \left( \int_0^T \left\| \sigma (\underline{s}, {_n X}_{\underline{s}}) \right\|^2 ds \right)^{2p} \right) \\
& \leq & C (x_0, p, T) \left( 1 + \int_0^T \left( \mathbb{E} | b (\underline{s}, {_n X}_{\underline{s}}) |^{4p} + \mathbb{E} \left\| \sigma (\underline{s}, {_n X}_{\underline{s}}) \right\|^{4p} \right) ds \right) \\
& \leq & C (x_0, p, T, K_{Lip}) \left( 1 + \int_0^T \mathbb{E} \| {_n X}_s \|_{\infty}^{4p} \: ds \right)
\end{eqnarray*}

We deduce from Gronwall's inequality that $ \mathbb{E} \| {_n X}_T \|_{\infty}^{4p} $ is bounded by a constant which depends only on $ x_0 $,  $ p $, $ T $ and $ K_{Lip} $.

Combining this result with \eqref{bm} and \eqref{cs}, we get:
\[ \mathbb{E} \| {_n \hat{X}}_T - {_n X}_T \|_{\infty}^{2p} \leq C (x_0, p, T, K_{Lip}) \left( \frac{1 + \log n}{n} \right)^p \]

Finally \eqref{sum} becomes:
\begin{eqnarray*}
&& \mathbb{E} \| X_T - {_n X}_T \|_{\infty}^{2p} \\
& \leq & C (p) \left( \mathbb{E} \| X_T - {_n \hat{X}}_T \|_{\infty}^{2p} + \mathbb{E} \| {_n \hat{X}}_T - {_n X}_T \|_{\infty}^{2p} \right) \\
& \leq & C (x_0, p, T, K_{Lip}) \left( \left( \frac{1 + \log n}{n} \right)^p + \int_0^T \mathbb{E} \| X_s - {_n X}_s \|_{\infty}^{2p} ds \right)
\end{eqnarray*}
And we conclude by Gronwall's inequality.
\qedhere
\end{proof}

\begin{cor} \label{cor_as1}
Under  Assumption \ref{hyp_sigma},
\[  \forall \alpha \in [0, \frac{1}{2}), \qquad n^{\alpha} \| X_T - {_n X}_T \|_{\infty} \underset{n \to \infty}{\longrightarrow} 0, \quad \mathbb{P} \mathrm{-a.s.} \]
\end{cor}
\begin{proof}
Let $ \alpha \in [0, \frac{1}{2}) $. For a $ p $ large enough, by \textbf{Proposition \ref{prop_euler}}, we have:
\[ \mathbb{E} \left[ \sum_{n \geq 1} n^{2 p \alpha} \| X_T - {_n X}_T \|_{\infty}^{2p} \right] < \infty \]
Thus
\[ \sum_{n \geq 1} n^{2 p \alpha} \| X_T - {_n X}_T \|_{\infty}^{2p} < \infty, \quad \mathbb{P} \mathrm{-a.s.} \]
and
\[ n^{\alpha} \| X_T - {_n X}_T \|_{\infty} \underset{n \to \infty}{\longrightarrow} 0, \quad \mathbb{P} \mathrm{-a.s.} \]
\qedhere
\end{proof}

%Gobet and Menozzi \cite{MR2409016} study weak convergence rates for Euler approximations for SDES with path-dependent coefficients.

\section{Smooth functional approximations for martingales} \label{sec:approximation}

Let $ g: D( [0, T], \mathbb{R}^d ) \longrightarrow \mathbb{R} $ be  a functional which satisfies:
\begin{hyp} \label{hyp_g1}
$ g: ( D( [0, T], \mathbb{R}^d ), \| \cdot \|_{\infty} ) \longrightarrow \mathbb{R} $ is continuous with  polynomial growth:
\[ \exists q \in \mathbb{N}, \exists C > 0, \forall \omega \in D( [0, T], \mathbb{R}^d ), \quad\left| g (\omega) \right| \leq C \left( 1 + \| \omega \|_{\infty}^q \right). \]
\end{hyp}
The (square-integrable) martingale:
\[ Y(t) = \mathbb{E} \left[ g ( X_T ) | \mathcal{F}_t \right] \]
may be represented as a non-anticipative functional of  $ W $:
\[ Y (t) = F (t, W_t) \]
where the functional $F$ is square-integrable but fail to be smooth a priori (see \cite{contriga2014} for results on pathwise differentiability of $F$). By \textbf{Proposition \ref{prop_mart}} we have:
\begin{eqnarray*}\label{eq:MRformula}
g ( X_T ) = Y(T) = Y(t) + \int_t^T \nabla_W Y (s) \cdot dW(s) \quad \mathbb{P} \mathrm{-a.s.}
\end{eqnarray*}
where $ \nabla_W Y $ is the weak vertical derivative of $ Y $ with respect to $ W $, which cannot be computed directly as a pathwise directional derivative unless $ F $ is a smooth functional (for example $ \in \mathbb{C}_{loc, r}^{1, 2} (\Lambda_T) $).   

The main idea is to approximate the martingale $ Y $ by a sequence of {\it smooth} martingales $ {_n Y}  $ which admit a functional representation $ {_n Y} (s) = F_n (s, W_s ) $ with $ F_n \in \mathbb{C}_{loc, r}^{1, 2} (\Lambda_T) $. Then by the functional It\^o formula, we have:
\[ \int_t^T \nabla_{\omega} F_n (s, W_s) \cdot dW(s) = {_n Y} (T) - {_n Y} (t) \underset{n \to \infty}{\longrightarrow} Y (T) - Y (t) = \int_t^T \nabla_W Y (s) \cdot dW (s). \]

One can then use the following estimator for $ \nabla_W Y $:
\[ Z_n (s) =  \nabla_{\omega} F_n (s, W_s), \]
where the vertical derivative $\nabla_{\omega} F_n (s, W_s) = \left( \partial_i F_n (s, W_s), 1 \leq i \leq d \right)$ may be computed as a  pathwise derivative
\[ \partial_i F_n (s, W_s) =  \lim\limits_{h \to 0} \frac{F_n (s, W_s + h e_i {\mathbf 1}_{[s, T]}) - F_n (s, W_s)}{h}, \]
yielding a concrete procedure for computing the estimator.

We will show in this section that the familiar {\it weak Euler approximation}  provides a systematic way of constructing such smooth functional approximations in the sense of Definition \ref{def_loc2}.

Define the concatenation of two c\`adl\`ag paths $ \omega, \omega' \in D ( [0,T], \mathbb{R}^d ) $ at time $ s \in [0, T] $, which we note $ \omega \conc{t} \omega' $, as the following c\`adl\`ag path on $ [0, T] $:
\[ \omega \conc{s} \omega' = \omega_s \conc{s} \omega' =
\begin{cases}
\omega (u) & u \in [0, s) \\
\omega (s) + \omega' (u) - \omega' (s) & u \in [s, T]
\end{cases}
\]
Observe that:
\[ \forall z \in \mathbb{R}^d,\quad \omega_s^z \conc{s} \omega' = ( \omega_s \conc{s} \omega' ) + z {\mathbf 1}_{[s, T]}.  \]

\begin{defi}[Weak Euler approximation]\label{def:weakEuler}
We define the (level-$n$) weak Euler approximation  of $F$ as the functional $F_n$ defined by
\begin{equation} \label{Fn}
F_n (s, \omega_s) = \mathbb{E} \left[ g ( {_n X} ( \omega_s \conc{s} W_T ) ) \right]
\end{equation}
\end{defi}
Applying this functional to the path of the Wiener process $W$, we obtain a $ (\mathcal{F}_t)_{t\geq 0}$-adapted process:\[ {_n Y} (s) =F_n (s, W_s). \]
Using independence of increments of $W$, we have
\[ {_n Y} (s) = \mathbb{E} \left[ g ( {_n X} ( W_T ) ) | \mathcal{F}_s \right] = \mathbb{E} \left[ g ( {_n X} ( W_s \conc{s} W_T ) ) | \mathcal{F}_s \right] = \mathbb{E} \left[ g ( {_n X} ( W_s \conc{s} B_T ) ) | \mathcal{F}_s \right] \]
where $B$ is any Wiener process independent from $W.$
In particular $ {_n Y}$ is a square-integrable martingale, so is weakly differentiable in the sense of \cite[Theorem 5.8]{CF13}. We will now show that  $F_n$ is in fact a smooth functional in the sense of Definition \ref{def_loc2}.
\begin{theorem} \label{th_diff}
Under Assumptions \ref{hyp_sigma} and  \ref{hyp_g1}, the functional $ F_n $ defined in \eqref{Fn} is horizontally differentiable and infinitely vertically differentiable.
\end{theorem}

\begin{proof}
First, note that under Assumption \ref{hyp_sigma}, ${_n X} (\omega)$ is bounded by a polynomial in the variables
$\omega (t_1-) - \omega (0) , \omega (t_2-) - \omega (t_{1}-),\cdots, \omega (t_n-) - \omega (t_{n-1}-))$. Combined with Assumption \ref{hyp_g1}, this implies that all expectations in the proof of this theorem are all well defined.

Let $ (s, \omega) \in \Lambda_T $ with $ t_k \leq s < t_{k+1} $ for some $ 0 \leq k \leq n-1 $. We start with the vertical differentiability of $ F_n $ at $ (s, \omega) $, which is equivalent to the differentiability at $ 0 $ of the  map $v:\mathbb{R}^d\to\mathbb{R}$ defined by:
\[ v (z) = F_n (s, \omega_s^z) = \mathbb{E} \left[ g ( {_n X} ( \omega_s^z \conc{s} B_T ) ) \right].\]
The main idea of the proof is to absorb the dependence with respect to $ z $ in the Gaussian density function when taking the expectation, which then implies smoothness.

As we have already shown, $ {_n X} ( \omega_s^z \conc{s} B_T ) $ depends only on $ ( \omega_s^z \conc{s} B_T ) (t_1-) - ( \omega_s^z \conc{s} B_T ) (0) $, $ \cdots $, $ ( \omega_s^z \conc{s} B_T ) (t_n-) - ( \omega_s^z \conc{s} B_T ) (t_{n-1}-) $, which are all explicit using the definition of the concatenation. For $ j < k $, we have:
\[ ( \omega_s^z \conc{s} B_T ) (t_{j+1}-) - ( \omega_s^z \conc{s} B_T ) (t_j-) = \omega (t_{j+1}-) - \omega (t_j-) \]
In the case where $ j = k $, we have:
\begin{eqnarray*}
&& ( \omega_s^z \conc{s} B_T ) (t_{k+1}-) - ( \omega_s^z \conc{s} B_T ) (t_k-) \\
&=& B (t_{k+1}) - B (s) + \omega (s) + z - \omega (t_k-) \\
&=& B (t_{k+1}) - B (s) + z + \omega (s) - \omega (t_k-)
\end{eqnarray*}
And for $ j > k $, we have:
\begin{eqnarray*}
&& ( \omega_s^z \conc{s} B_T ) (t_{j+1}-) - ( \omega_s^z \conc{s} B_T ) (t_j-) \\
&=& B (t_{j+1}) - B (s) + \omega (s) + z - ( B (t_j) - B (s) + \omega (s) + z ) \\
&=& B (t_{j+1}) - B (t_j)
\end{eqnarray*}
Thus we have:
\begin{eqnarray*}
{_n X} ( \omega_s^z \conc{s} B_T ) =& p_n \Big( \omega (t_1-) - \omega (0), \cdots, \omega (t_k-) - \omega (t_{k-1}-), B (t_{k+1}) - B (s) + z + \omega (s) - \omega (t_k-), \\
& B (t_{k+2}) - B (t_{k+1}), \cdots, B (t_n) - B (t_{n-1}) \Big)
\end{eqnarray*}
where $ p_n: \mathcal{M}_{d, n} (\mathbb{R}) \rightarrow D ([0, T], \mathbb{R}^d) $ is the map defined by \eqref{p}.

Observe from the previous equation that, for a fixed $ z $, the value of $ {_n X} ( t_{k+1}, \omega_s^z \conc{s} B_T ) $ as a random variable depends only on a finite number of Gaussian variables: $ B (t_{k+1}) - B (s) $, $ B (t_{k+2}) - B (t_{k+1}) $, $ \cdots $, $ B (t_j) - B (t_{j-1}) $. Since the joint distribution of these Gaussian variables is explicit, $ v (z) = \mathbb{E} \left[ g ( {_n X}_T ( \omega_s^z \conc{s} B_T ) ) \right] $ can be computed explicitly as an integral in finite dimension.

Denoting $ y = (y_1, \cdots, y_{n-k}) \in \mathcal{M}_{d, n-k} ( \mathbb{R} ) $,
\begin{eqnarray}
\nonumber v(z) &=& \mathbb{E} \left[ g ( {_n X}_T ( \omega_s^z \conc{s} B_T ) ) \right] \\
\nonumber &=& \mathbb{E} \Big[ g \Big( p_n (\omega (t_1-) - \omega (0), \cdots, \omega (t_k-) - \omega (t_{k-1}-), \\
\nonumber && B (t_{k+1}) - B (s) + z + \omega (s) - \omega (t_k-), B (t_{k+2}) - B (t_{k+1}), \cdots, B (t_n) - B (t_{n-1}) \Big) \Big] \\
\nonumber &=& \int_{\mathbb{R}^{d \times (n-k)}} g \Big ( p_n (\omega (t_1-) - \omega (0), \cdots, \omega (t_k-) - \omega (t_{k-1}-), y_1 + z + \omega (s) - \omega (t_k-), \\
\nonumber && y_2, \cdots, y_{n-k}) \Big) \Phi (y_1, t_{k+1} - s) \prod_{l=2}^{n-k} \Phi (y_l, \delta) dy_1 dy_2 \cdots dy_{n-k} \\
\nonumber &=& \int_{\mathbb{R}^{d \times (n-k)}} g \Big ( p_n (\omega (t_1-) - \omega (0), \cdots, \omega (t_k-) - \omega (t_{k-1}-), y_1 + \omega (s) - \omega (t_k-), \\
\label{vertical} && y_2, \cdots, y_{n-k}) \Big) \Phi (y_1 - z, t_{k+1} - s) \prod_{l=2}^{n-k} \Phi (y_l, \delta) dy_1 dy_2 \cdots dy_{n-k}
\end{eqnarray}
with
\[ \Phi (x, t) = ( 2 \pi t )^{-\frac{d}{2}} \exp \left( -\frac{|x|^2}{2 t} \right), \quad x \in \mathbb{R}^d \]
the  density function of a $ d $-dimensional $N(0,t I_d)$ variable.
Since the only term which depends on $ z $ in the integrand of  \eqref{vertical} is $ \Phi (y_1 - z, t_{k+1} - s) $, which is a smooth function of $ z $, thus $ v $ is differentiable at all $ z \in \mathbb{R}^d $, in particular at $ 0 $. Hence $ F_n $ is vertically differentiable at $ (s, \omega) \in \Lambda_T $ with: for $ 1 \leq i \leq d $,
\begin{eqnarray}
\nonumber \partial_i F_n (s, \omega) &=& \int_{\mathbb{R}^{d \times (n-k)}} g \Big ( p_n (\omega (t_1-) - \omega (0), \cdots, \omega (t_k-) - \omega (t_{k-1}-), y_1 + \omega (s) - \omega (t_k-), \\
\nonumber && y_2, \cdots, y_{n-k}) \Big) \frac{y_1 \cdot e_i}{t_{k+1} - s} \Phi(y_1, t_{k+1} - s) \prod_{l=2}^{n-k} \Phi (y_l, \delta) dy_1 dy_2 \cdots dy_{n-k} \\
\label{MC} &=& \mathbb{E} \left[ g ( {_n X} (\omega_s \conc{s} B_T) ) \frac{( B (t_{k+1}) - B (s) ) \cdot e_i}{t_{k+1} - s} \right]
\end{eqnarray}

Remark that when $ s $ tends towards $ t_{k+1} $, $ \nabla_{\omega} F_n (s, \omega) $ may tend to infinity because of the term $ t_{k+1} - s $ in the denominator. However in the interval $ [t_k, t_{k+1}) $, $ \nabla_{\omega} F_n (s, \omega) $ behaves well and is locally bounded.

Iterating this procedure, one can show that $ F_n $ is  vertically differentiable to any order. For example, we have: for $ z \in \mathbb{R}^d $,
\begin{eqnarray*}
\partial_i F_n (s, \omega_s^z) &=& \int_{\mathbb{R}^{d \times (n-k)}} g \Big ( p_n (\omega (t_1-) - \omega (0), \cdots, \omega (t_k-) - \omega (t_{k-1}-), y_1 + z + \omega (s) - \omega (t_k-), \\
&& y_2, \cdots, y_{n-k}) \Big) \frac{y_1 \cdot e_i}{t_{k+1} - s} \Phi(y_1, t_{k+1} - s) \prod_{l=2}^{n-k} \Phi (y_l, \delta) dy_1 dy_2 \cdots dy_{n-k}
\end{eqnarray*}
Thus we have:
\begin{eqnarray*}
\partial_i^2 F_n (s, \omega) &=& \int_{\mathbb{R}^{d \times (n-k)}} g \Big( p_n (\omega (t_1-) - \omega (0), \cdots, \omega (t_k-) - \omega (t_{k-1}-), y_1 + \omega (s) - \omega (t_k-), y_2, \cdots, y_{n-k}) \Big) \\
&& \left( \frac{(y_1 \cdot e_i)^2}{(t_{k+1}-s)^2} - \frac{1}{t_{k+1} - s} \right) \Phi(y_1, t_{k+1} - s) \prod_{l=2}^{n-k} \Phi (y_l, \delta) dy_1 dy_2 \cdots dy_{n-k}
\end{eqnarray*}
And for $ i \neq j $:
\begin{eqnarray*}
\partial_{ij} F_n (s, \omega) &=& \int_{\mathbb{R}^{d \times (n-k)}} g \Big( p_n (\omega (t_1-) - \omega (0), \cdots, \omega (t_k-) - \omega (t_{k-1}-), y_1 + \omega (s) - \omega (t_k-), \\
&& y_2, \cdots, y_{n-k}) \Big) \frac{(y_1 \cdot e_i) (y_1 \cdot e_j)}{(t_{k+1}-s)^2} \Phi(y_1, t_{k+1} - s) \prod_{l=2}^{n-k} \Phi (y_l, \delta) dy_1 \cdots dy_{n-k}
\end{eqnarray*}

The horizontal differentiability of $ F_n $ can be proved similarly. Consider the following map:
\[ w (h) = F_n (s+h, \omega_s) = \mathbb{E} \left[ g ( {_n X} ( \omega_s \conc{s+h} B_T ) ) \right], \quad h > 0 \]
The objective is to show that $ w $ is right-differentiable at $ 0$.

We assume again that $ t_k \leq s < t_{k+1} $ for some $ 0 \leq k \leq n-1 $, and we take an $ h > 0 $ small enough such that $ s + h < t_{k+1} $. Using the same argument and the fact that $ \omega_s (s+h) = \omega (s) $, we have:
\begin{eqnarray*}
{_n X} ( \omega_s \conc{s+h} B_T ) &=& p_n (\omega (t_1-) - \omega (0), \cdots, \omega (t_k-) - \omega (t_{k-1}-), B (t_{k+1}) - B (s+h) + \omega (s) - \omega (t_k-), \\
&& B (t_{k+2}) - B (t_{k+1}), \cdots, B (t_n) - B (t_{n-1}))
\end{eqnarray*}

Let $ y = (y_1, \cdots, y_{n-k}) \in \mathcal{M}_{d, n-k} ( \mathbb{R} ) $. We calculate explicitly $ w (h) $:
\begin{eqnarray}
\nonumber w(h) &=& \mathbb{E} \left[ g ( {_n X}_T ( \omega_s \conc{s+h} B_T ) ) \right] \\
\nonumber &=& \mathbb{E} \Big[ g \Big( p_n (\omega (t_1-) - \omega (0), \cdots, \omega (t_k-) - \omega (t_{k-1}-), B (t_{k+1}) - B (s+h) + \omega (s) - \omega (t_k-), \\
\nonumber && B (t_{k+2}) - B (t_{k+1}), \cdots, B (t_n) - B (t_{n-1})) \Big) \Big] \\
\nonumber &=& \int_{\mathbb{R}^{d \times (n-k)}} g \Big( p_n (\omega (t_1-) - \omega (0), \cdots, \omega (t_k-) - \omega (t_{k-1}-), y_1 + \omega (s) - \omega (t_k-), \\
\label{horizontal} && y_2, \cdots, y_n) \Big) \Phi(y_1, t_{k+1}-s-h) \prod_{l=2}^{n-k} \Phi(y_l, \delta) dy_1 dy_2 \cdots dy_{n-k}
\end{eqnarray}

Again the only term which depends on $ h $ in the integrand of \eqref{horizontal} is $ \Phi (y_1, t_{k+1}-s-h) $, which is a smooth function of $ h $. Therefore $ F_n $ is horizontally differentiable with:
\begin{eqnarray*}
\mathcal{D} F_n (s, \omega_s) &=& \int_{\mathbb{R}^{d \times (n-k)}} g \Big( p_n(\omega (t_1-) - \omega (0), \cdots, \omega (t_k-) - \omega (t_{k-1}-), y_1 + \omega (s) - \omega (t_k-), \\
&& y_2, \cdots, y_n) \Big) \left( \frac{d}{2 (t_{k+1}-s)} - \frac{|y_1|^2}{2 (t_{k+1}-s)^2} \right) \Phi(y_1, t_{k+1}-s) \prod_{l=2}^{n-k} \Phi(y_l, \delta) dy_1 \cdots dy_{n-k}.
\end{eqnarray*}
\qedhere
\end{proof}

The following result shows that the functional derivatives of  $ F_n $ satisfy the necessary regularity conditions for applying the functional It\^o formula to $ F_n $:
\begin{theorem} \label{th_reg}
Under Assumptions \ref{hyp_sigma} and  \ref{hyp_g1}, $ F_n \in \mathbb{C}^{1, 2}_{loc, r} (\Lambda_T) $.
\end{theorem}

\begin{proof}
We have already shown in Theorem \ref{th_diff} that $ F_n $ is horizontally differentiable and twice vertically differentiable. Using the expressions of $ \mathcal{D} F_n $, $ \nabla_{\omega} F_n $ and $ \nabla^2_{\omega} F_n $ obtained in the proof of \ref{th_diff} and the assumption that $ g $ has at most polynomial growth with respect to $ \| \cdot \|_{\infty} $,  we observe that in each interval $ [t_k, t_{k+1}) $ with $ 0 \leq k \leq n-1 $, $ \mathcal{D} F_n $, $ \nabla_{\omega} F_n $ and $ \nabla^2_{\omega} F_n $ satisfy the boundedness-preserving property. We now prove that $ F_n $ is left-continuous, $ \nabla_{\omega} F_n $ and $ \nabla^2_{\omega} F_n $ are right-continuous, and $ \mathcal{D} F_n $ is continuous at fixed times.

Let $ s \in [t_k, t_{k+1}) $ for some $ 0 \leq k \leq n-1 $ and $ \omega \in D ([0, T], \mathbb{R}^d) $. We first prove that $ F_n $ is right-continuous at $ (s, \omega) $, and is jointly continuous at $ (s, \omega) $ for $ s \in (t_k, t_{k+1}) $. By definition of joint-continuity (or right-continuous), we want to show that: $ \forall \epsilon > 0, \exists \eta > 0, \forall (s', \omega') \in \Lambda_T $ (for the right-continuity, we assume in addition that $ s' > s $),
\[ d ( (s, \omega), (s', \omega') ) < \eta ) \Rightarrow | F_n (s, \omega) - F_n (s', \omega') | < \epsilon \]

Let $ (s', \omega') \in \Lambda_T $ (with $ s' > s $ for the right-continuity). We assume that $ d ( (s, \omega), (s', \omega') ) \leq \eta $ with an $ \eta $ small enough such that $ s' \in [t_k, t_{k+1}) $ (this is always possible as if $ s = t_k $, we are only interested in the right-continuity, thus $ s' > s $). It suffices to prove that $ | F_n (s, \omega) - F_n (s', \omega') | \leq C (s, \omega_s, \eta) $ with $ C (s, \omega_s, \eta) $ a quantity depending only on $ s $, $ \omega_s $ and $ \eta $, and $ C (s, \omega_s, \eta) \underset{\eta \to 0}{\longrightarrow} 0 $.

We use the expression of $ F_n $ obtained in the proof of  {Theorem \ref{th_diff}}. Denote $ y = (y_1, \cdots, y_{n-k}) \in \mathcal{M}_{d, n-k} ( \mathbb{R} ) $. We have
\begin{eqnarray*}
F_n (s, \omega) &=& \int_{\mathbb{R}^{d \times (n-k)}} g \Big( p_n (\omega (t_1-) - \omega (0), \cdots, \omega (t_k-) - \omega (t_{k-1}-), y_1 + \omega (s) - \omega (t_k-), \\
&& y_2, \cdots, y_n) \Big) \Phi(y_1, t_{k+1}-s) \prod_{l=2}^{n-k} \Phi(y_l, \delta) dy_1 dy_2 \cdots dy_{n-k}
\end{eqnarray*}
and
\begin{eqnarray*}
F_n (s', \omega') &=& \int_{\mathbb{R}^{d \times (n-k)}} g \Big( p_n (\omega' (t_1-) - \omega' (0), \cdots, \omega' (t_k-) - \omega' (t_{k-1}-), y_1 + \omega' (s') - \omega' (t_k-), \\
&& y_2, \cdots, y_n) \Big) \Phi(y_1, t_{k+1}-s') \prod_{l=2}^{n-k} \Phi(y_l, \delta) dy_1 dy_2 \cdots dy_{n-k}.
\end{eqnarray*}

To simplify the notations, we set:
\[ \tilde{p} (\omega, s, y) = p_n (\omega (t_1-) - \omega (0), \cdots, \omega (t_k-) - \omega (t_{k-1}-), y_1 + \omega (s) - \omega (t_k-), y_2, \cdots, y_n) \]
and
\[ \tilde{p} (\omega', s', y) = p_n (\omega' (t_1-) - \omega' (0), \cdots, \omega' (t_k-) - \omega' (t_{k-1}-), y_1 + \omega' (s') - \omega' (t_k-), y_2, \cdots, y_n) \]
Similarly $ \tilde{p}_t ( \cdot ) $ will be the path of $ \tilde{p} ( \cdot ) $ stopped at time $ t $.

As $ \| \omega_s - \omega_{s'} \|_{\infty} \leq \eta $, by \textbf{Lemma \ref{lem_euler}}, we have:
\[ \| \tilde{p} (\omega, s, y) - \tilde{p} (\omega', s', y) \|_{\infty} \leq C (\omega_s, y, K_{Lip}, T) \eta \]

We now control the difference between $ F_n (s, \omega) $ and $ F_n (s', \omega') $.
\begin{eqnarray}
\nonumber && | F_n (s, \omega) - F_n (s', \omega') | \\
\nonumber & \leq & \int_{\mathbb{R}^{d \times (n-k)}} | g ( \tilde{p} (\omega, s, y) ) \Phi(y_1, t_{k+1}-s) - g ( \tilde{p} (\omega', s', y) ) \Phi(y_1, t_{k+1}-s') | \prod_{l=2}^{n-k} \Phi(y_l, \delta) dy \\
\nonumber & \leq & \int_{\mathbb{R}^{d \times (n-k)}} \Big( | g ( \tilde{p} (\omega, s, y) ) - g ( \tilde{p} (\omega', s', y) ) | \Phi(y_1, t_{k+1}-s') \\
\label{control} &&+ | g ( \tilde{p} (\omega, s, y) ) | \cdot | \Phi(y_1, t_{k+1}-s) - \Phi(y_1, t_{k+1}-s') | \Big) \times \prod_{l=2}^{n-k} \Phi(y_l, \delta) dy_1 dy_2 \cdots dy_{n-k}
\end{eqnarray}

Observe that $ | \Phi(y_1, t_{k+1}-s) - \Phi(y_1, t_{k+1}-s') | \leq |s-s'| \cdot \rho (y_1, \eta) \leq \rho (y_1, \eta) \cdot \eta^2 $ with
\[ \rho (y_1, \eta) = \sup_{t \in [t_{k+1}-s - \eta^2, \delta]} | \partial_t \Phi(y_1, t) | \]
and we have:
\[ \rho (y_1, \eta) \underset{\eta \to 0}{\longrightarrow} \sup_{t \in [t_{k+1}-s, \delta]} | \partial_t \Phi(y_1, t) | = \sup_{t \in [t_{k+1}-s, \delta]} \left| \Phi(y_1, t) \left( \frac{|y_1|^2}{2 t^2} - \frac{d}{2 t} \right) \right| < \infty \]
So the second part of \eqref{control} can be controlled by:
\[ \int_{\mathbb{R}^{d \times (n-k)}} | g ( \tilde{p} (\omega, s, y) ) | \cdot | \Phi(y_1, t_{k+1}-s) - \Phi(y_1, t_{k+1}-s') | \prod_{l=2}^{n-k} \Phi(y_l, \delta) dy_1 dy_2 \cdots dy_{n-k} \leq C (s, \omega_s, \eta) \]
with
\[ C (s, \omega_s, \eta) \underset{\eta \to 0}{\longrightarrow} 0. \]
For the first part of \eqref{control}, we use the continuity of $ g $. As $ g $ is continuous at $ p (\omega_s, y) $, we have:
\[ | g ( \tilde{p} (\omega, s, y) ) - g ( \tilde{p} (\omega', s', y) ) | \leq C (s, \omega_s, y, \eta) \]
with
\[ C (s, \omega_s, y, \eta) \underset{\eta \to 0}{\longrightarrow} 0 \]
and
\[ \Phi(y_1, t_{k+1}-s') \leq \sup_{t \in [t_{k+1}-s-\eta^2, \delta]} \Phi(y_1,t) < \infty. \]
Thus the first part of \eqref{control} can also be bounded by $ C (s, \omega_s, \eta) $ with
\[ C (s, \omega_s, \eta) \underset{\eta \to 0}{\longrightarrow} 0. \]
We conclude that $ | F_n (s, \omega_s) - F_n (s', \omega'_{s'}) | \leq C (s, \omega_s, \eta) $ depending only on $ s $, $ \omega_s $ and $ \eta $, and $ C (s, \omega_s, \eta) \underset{\eta \to 0}{\longrightarrow} 0 $, which proves the right-continuity of $ F_n $ and the joint-continuity of $ F_n $ at all $ (s, \omega) \in \Lambda_T $ for $ s \neq t_k $, $ 0 \leq k \leq n-1 $.

The right-continuity of $ \nabla_{\omega} F_n $, $ \nabla_{\omega}^2 F_n $  and the continuity at fixed times of $ \mathcal{D} F_n $ can be similarly shown from the expressions of $\nabla_{\omega} F_n $, $ \nabla_{\omega}^2 F_n $ and $ \mathcal{D} F_n $ obtained in  Theorem \ref{th_diff}.
Now it remains to show that for $ 1 \leq k \leq n $ and $ \omega \in D ([0, T], \mathbb{R}^d) $, $ F_n $ is left-continuous at $ (t_k, \omega) $. Let $ (s', \omega') \in \Lambda_T $ with $ s' < t_k $ such that $ d ((t_k, \omega), (s', \omega')) \leq \eta $. We choose an $  \eta $ small enough in order that $ s' \in [t_{k-1}, t_k) $, and we want to show that $ | F_n (t_k, \omega) - F_n (s', \omega') | \leq C (t_k, \omega_{t_k}, \eta) $ for some $ C (t_k, \omega_{t_k}, \eta) $ depending only on $ t_k $, $ \omega_{t_k} $ and $ \eta $ with $ C (t_k, \omega_{t_k}, \eta) \underset{\eta \to 0}{\longrightarrow} 0 $.

We first decompose $ | F_n (t_k, \omega) - F_n (s', \omega') | $ into two terms:
\[ | F_n (t_k, \omega) - F_n (s', \omega') | \leq | F_n (t_k, \omega) - F_n  (s', \omega_{s'}) | + | F_n  (s', \omega_{s'}) - F_n (s', \omega') | \]
For the second part, as $ F_n $ is continuous at fixed time $ s' $ by the first part of the proof, and $ \|  \omega_{s'} - \omega'_{s'} \|_{\infty} \leq \eta $, we have $ | F_n  (s', \omega_{s'}) - F_n (s', \omega') | \leq C (t_k, \omega_{t_k}, \eta) $ with $ C (t_k, \omega_{t_k}, \eta) \underset{\eta \to 0}{\longrightarrow} 0 $.

For the first part $ | F_n (t_k, \omega) - F_n (s', \omega_{s'}) | $, the difficulty is that $ s' $ and $ t_k $ no longer lie in the same interval, thus we need to perform one more integration for $ F_n (s', \omega_{s'}) $ compared to $ F_n (t_k, \omega) $. Let $ y = (y_1, \cdots, y_{n-k}) \in \mathcal{M}_{d, n-k} (\mathbb{R}) $ and $ y' \in \mathbb{R}^d $. Using again the expression of $ F_n $ we have obtained in the proof of  {Theorem \ref{th_diff}}, we have:
\begin{eqnarray*}
F_n (t_k, \omega) &=& \int_{\mathbb{R}^{d \times (n-k)}} g \Big ( p_n (\omega (t_1-) - \omega (0), \cdots, \omega (t_k-) - \omega (t_{k-1}-), y_1 + \omega (t_k) - \omega (t_k-), \\
&& y_2, \cdots, y_{n-k}) \Big) \prod_{l=1}^{n-k} dy_l \Phi (y_l, \delta) 
\end{eqnarray*}
and
\begin{eqnarray*}
F_n (s', \omega_{s'}) &=& \int_{\mathbb{R}^{d \times (n-k+1)}} g \Big( p_n (\omega (t_1-) - \omega (0), \cdots, y' + \omega (s') - \omega (t_{k-1}-), y_1, \cdots, y_n) \Big) \\
&& \Phi (y', t_k - s') dy'  \prod_{l=1}^{n-k} dy_l \Phi (y_l, \delta) \\
&=& \int_{\mathbb{R}^d} \Big( \int_{\mathbb{R}^{d \times (n-k)}} g ( p (\omega (t_1-) - \omega (0), \cdots, y' + \omega (s') - \omega (t_{k-1}-), y_1, \cdots, y_n) ) \\
&& \prod_{l=1}^{n-k} \Phi (y_l, \delta) dy_l \Big) \Phi (y', t_k - s') dy'
\end{eqnarray*}

We now define $ \zeta: \mathbb{R}^d \rightarrow \mathbb{R} $ by: for $ y' \in \mathbb{R}^d $,
\[ \zeta (y') = \int_{\mathbb{R}^{d \times (n-k)}} g ( p_n (\omega (t_1-) - \omega (0), \cdots, y' + \omega (s') - \omega (t_{k-1}-), y_1, \cdots, y_n) ) \prod_{l=1}^{n-k} \Phi (y_l, \delta) dy_l \]

By \textbf{Lemma \ref{lem_euler}} and the continuity of $ g $ with respect to $ \| \cdot \|_{\infty} $, the map
\[ y' \mapsto g ( p (\omega (t_1-) - \omega (0), \cdots, y' + \omega (s') - \omega (t_{k-1}-), y_1, \cdots, y_n) ) \]
is continuous. As $ g $ has at most polynomial growth with respect to $ \| \cdot \|_{\infty} $, by the dominated convergence theorem, $ \zeta $ is also continuous. Moreover, $ \zeta $ has at most polynomial growth. And as $ t_k - s' \leq \eta^2 $, we have
\begin{eqnarray*}
F_n (s', \omega_{s'}) &=& \int_{\mathbb{R}^d} \zeta (y') \Phi (y', t_k - s') dy' \\
&=& \int_{\mathbb{R}^d} ( \zeta (y') - \zeta (0) ) \Phi (y', t_k - s') dy' + \zeta (0).
\end{eqnarray*}
with
\[ \left| \int_{\mathbb{R}^d} ( \zeta (y') - \zeta (0) ) \Phi (y', t_k - s') dy' \right| \leq C (t_k, \omega_{t_k}, \eta) \]
and $ C (t_k, \omega_{t_k}, \eta) \underset{\eta \to 0}{\longrightarrow} 0 $.

It remains to control the difference between $ F_n (t_k, \omega) $ and $ \zeta (0) $. We remark that:
\begin{eqnarray*}
\zeta (0) &=& \int_{\mathbb{R}^{d \times (n-k)}} g ( p_n (\omega (t_1-) - \omega (0), \cdots, \omega (s') - \omega (t_{k-1}-), y_1, \cdots, y_n) ) \prod_{l=1}^{n-k} \Phi (y_l, \delta) dy_l  \\
&=& \mathbb{E} \left[ g ({_n X}_T (\omega_{s'} \conc{t_k} B_T)) \right] = F_n (t_k, \omega_{s'})
\end{eqnarray*}
As $ \| \omega_{s'} - \omega_{t_k} \|_{\infty} \leq \| \omega_{s'} - \omega'_{s'} \|_{\infty} + \| \omega_{t_k} - \omega'_{s'} \|_{\infty} \leq 2 \eta $, again by the continuity of $ F_n $ at fixed time $ t_k $ established in the first part of the proof, we have:
\[ | F_n (t_k, \omega) - \zeta (0) | \leq C (t_k, \omega_{t_k}, \eta) \]
with $ C (t_k, \omega_{t_k}, \eta) \underset{\eta \to 0}{\longrightarrow} 0 $.

We conclude that
\[ | F_n (t_k, \omega) - F_n (s', \omega') | \leq C (t_k, \omega_{t_k}, \eta) \]
with $ C (t_k, \omega_{t_k}, \eta) \underset{\eta \to 0}{\longrightarrow} 0 $, which proves the left-continuity of $ F_n $ at $ (t_k, \omega) $.
\qedhere
\end{proof}

\begin{cor} \label{cor_ito}
Under Assumptions \ref{hyp_sigma} and   \ref{hyp_g1}, for any $ t \in [0, T) $ we have:
\begin{equation} \label{ito_c}
F_n (T, W_T) - F_n (t, W_t) = \int_t^T \nabla_{\omega} F_n (s, W_s) \cdot dW (s), \quad \mathbb{P}-a.s.
\end{equation}
\end{cor}
\begin{proof}
As $ F_n \in \mathbb{C}^{1, 2}_{loc, r} (\Lambda_T) $, we can apply the functional It\^o formula \textbf{Proposition \ref{prop_ito2}} and we remark that the finite variation term is zero as $ {_n Y} (s) = F_n (s, W_s) $ is a martingale.
\qedhere
\end{proof}

\begin{rmq} \label{rmq_mart}
We can also verify using directly the expressions we have obtained in  {Theorem \ref{th_diff}} for $ \mathcal{D} F_n $ and $ \nabla_{\omega}^2 F_n $ that the finite variation terms in \eqref{ito_c}  cancel each other. By the functional Itô formula, the finite variation term in \eqref{ito_c} equals to $ \mathcal{D} F_n (s, W_s) + \frac{1}{2} \mathrm{tr} ( \nabla^2_{\omega} F_n (s, W_s) ) $. And for $ (s, \omega) \in \Lambda_T $ with $ s \in [t_k, t_{k+1}) $, we have:
\begin{eqnarray*}
&& \mathrm{tr} \left( \nabla_{\omega}^2 F_n (s, \omega) \right) = \sum_{i=1}^d \partial_i^2 F_n (s, \omega) \\
&=& \int_{\mathbb{R}^{d \times (n-k)}} g ( p_n (\omega (t_1-) - \omega (0), \cdots, \omega (t_k-) - \omega (t_{k-1}-), y_1 + \omega (s) - \omega (t_k-), y_2, \cdots, y_n) ) \\
&& \sum_{i=1}^d \left( \frac{(y_1 \cdot e_i)^2}{(t_{k+1}-s)^2} - \frac{1}{t_{k+1} - s} \right) \Phi(y_1, t_{k+1} - s) \prod_{l=2}^{n-k} \Phi (y_l, \delta) dy_1 \cdots dy_{n-k} \\
&=& \int_{\mathbb{R}^{d \times (n-k)}} g ( p_n (\omega (t_1-) - \omega (0), \cdots, \omega (t_k-) - \omega (t_{k-1}-), y_1 + \omega (s) - \omega (t_k-), y_2, \cdots, y_n) ) \\
&& \left( \frac{|y_1|^2}{(t_{k+1}-s)^2} - \frac{d}{t_{k+1}-s} \right) \Phi(y_1, t_{k+1} - s) \prod_{l=2}^{n-k} \Phi (y_l, \delta) dy_1 \cdots dy_{n-k} \\
&=& -2 \mathcal{D} F_n (s, \omega_s)
\end{eqnarray*}
which confirms that $F_n$ is a solution of the path-dependent Kolmogorov equation \cite[Sec. 5]{Cont2012}:
$$ \mathcal{D} F_n (s, W_s) + \frac{1}{2} \mathrm{tr} ( \nabla^2_{\omega} F_n (s, W_s) ) = 0.$$
\end{rmq}
\section{Convergence and error analysis}\label{sec:convergence}

In this section, we analyze the convergence rate of our approximation method. After having constructed a sequence of smooth functionals $ F_n $ ( {Theorem \ref{th_diff}} and  {Theorem \ref{th_reg}}), we can now approximate $ \nabla_W Y $ by:
\[ Z_n (s) =  \nabla_{\omega} F_n (s, W_s) \]
which, in contrast to the weak derivative $ \nabla_W Y $, is computable as  a pathwise directional derivative. In practice, $ \nabla_{\omega} F_n (s, W_s) $ can be computed numerically via a finite difference method or a Monte-Carlo method using the expression \eqref{MC} of $ \nabla_{\omega} F_n $.

For $ t \in [0, T] $, the quantity of interest is the integral of $ \nabla_W Y - Z_n $ along the path of $ W $ between $ t $ and $ T $, i.e.
\[ \int_t^T \left( \nabla_W Y - Z_n \right) \cdot dW = \int_t^T \nabla_W Y (s) \cdot dW (s) - \int_t^T \nabla_{\omega} F_n (s, W_s) \cdot dW (s) \]
By the martingale representation formula \textbf{Proposition \ref{prop_mart}} and \textbf{Corollary \ref{cor_ito}}, we have $ \mathbb{P} $-a.s.
\begin{eqnarray*}
 \int_t^T \left( \nabla_W Y - Z_n \right) \cdot dW 
&=& Y (T) - Y(t) - ({_n Y} (T) - {_n Y} (t)) \\
&=& g (X_T) - g ( {_n X}_T (W_T) ) - \mathbb{E} \left[ g(X_T) - g ( {_n X}_T (W_t \conc{t} B_T) ) | \mathcal{F}_t \right]
\end{eqnarray*}
where $ {_n X} $ is the path of the piecewise constant Euler-Maruyama scheme defined in \eqref{euler}.
Remark that by definition of the concatenation operation, if $ B $ and $ W $ are two independent Brownian motions, we have:
\begin{eqnarray*}
 \mathbb{E} \left[ g ( {_n X}_T (W_t \conc{t} B_T) ) | \mathcal{F}_t \right] 
&=& \mathbb{E} \left[ g ( {_n X}_T (W_t \conc{t} W_T) ) | \mathcal{F}_t \right] 
= \mathbb{E} \left[ g \left( {_n X}_T (W_T)\right)  | \mathcal{F}_t \right]
\end{eqnarray*}
\begin{prop} \label{cor_as2}
Under Assumptions \ref{hyp_sigma} and \ref{hyp_g1},
\[ \forall t \in [0, T], \quad \int_t^T \left( \nabla_W Y - Z_n \right) \cdot dW \underset{n \to \infty}{\longrightarrow} 0, \quad \mathbb{P} \mathrm{-a.s.} \]
\end{prop}
\begin{proof}
We have already shown that:
\[ \int_t^T \left( \nabla_W Y - Z_n \right) \cdot dX = g (X_T) - g ({_n X}_T) - \mathbb{E} \left[ g (X_T) - g ({_n X}_T) | \mathcal{F}_t \right] \]
As $ g $ is continuous with respect to $ \| \cdot \|_{\infty} $, by \textbf{Corollary \ref{cor_as1}}, we have:
\[ g (X_T) - g ({_n X}_T) \underset{n \to \infty}{\longrightarrow} 0, \quad \mathbb{P} \mathrm{-a.s.} \] 
Under assumption  \ref{hyp_g1}, the sequence of martingales $Y^n=E[g ({_n X}_T) |{\cal F_t}$ is bounded in $L^2$ therefore uniformly integrable, hence the result.
\qedhere
\end{proof}

\begin{cor}\label{cor_norm}
Under  Assumptions \ref{hyp_sigma} and   \ref{hyp_g1},
\[ \forall t\in [0,T],\quad \left\| \int_t^T \left( \nabla_W Y - Z_n \right) \cdot dW \right\|_{2p} \underset{n \to \infty}{\longrightarrow} 0, \quad \forall p \geq 1 \]
\end{cor}
\begin{proof}  $ g $ has at most polynomial growth with respect to $ \| \cdot \|_{\infty} $, which ensures the uniform integrability of $ | g ({_n X}_T) |^{2p} $. We conclude by applying the dominated convergence theorem. 
\end{proof}

Under a slightly stronger assumption on $ g $ we can obtain a rate of convergence:
\begin{theorem}[Rate of convergence] \label{th_error}
Let $ p \geq 1$ and assume  $ g: (D ([0, T], \mathbb{R}^d),\|.\|_\infty)\to\mathbb{R}$ is Lipschitz-continuous:
\[ \exists  g_{Lip} > 0,\quad \forall \omega, \omega' \in D ([0, T], \mathbb{R}^d),\quad \left| g (\omega) - g (\omega') \right| \leq g_{Lip} \:  \| \omega - \omega' \|_\infty. \]
Under  Assumptions \ref{hyp_sigma}, the $L^{2p} $-error of the approximation $  Z_n $ of $ \nabla_W Y $ along the path of $ W $ between
$ t $ and $ T $ is bounded by:
\[  \mathbb{E} \left( \left\| \int_t^T \left( \nabla_W Y - Z_n \right) \cdot dW \right\|^{2p}\right) \leq C (x_0, p, T, K_{Lip}, g_{Lip}) \left(\frac{1+\log n}{n}\right)^p, \quad \forall p \geq 1 \]
where the constant $ C $  depends only on $ x_0, p, T, K_{Lip}$ and $ g_{Lip} $. In particular:
\[ \forall \alpha \in [0, \frac{1}{2}),\qquad n^{\alpha} \left( \int_t^T (\nabla_W Y  - Z_n) \cdot dW \right) \underset{n \to \infty}{\longrightarrow} 0, \quad \mathbb{P} \mathrm{-a.s.} \]
\end{theorem}
\begin{proof}
This result is a consequence of Proposition \ref{prop_euler}  since
\begin{eqnarray*}
\left\| \int_t^T \left( \nabla_W Y (s) - Z_n (s) \right) \cdot dW (s) \right\|_{2p}
& \leq & \| g (X_T) - g ( {_n X}_T ) \|_{2p} + \| \mathbb{E} [ g(X_T) - g ( {_n X}_T ) | \mathcal{F}_t ] \|_{2p} \\
& \leq & 2 \| g (X_T) - g ( {_n X}_T ) \|_{2p} \\
& \leq & 2 g_{Lip} \| \sup_{s \in [0, T]} | X (s) - {_n X} (s) | \|_{2p}.
\end{eqnarray*}
\qedhere
\end{proof}

The following example how our result may be used to construct explicit approximations with controlled convergence rates for conditional expectation of non-smooth functionals:
\begin{example}
Let  $$ g(\omega) = \psi ( \omega(T),\mathop{\sup}_{t\in [0,T]} \|\omega(t)\| ) $$ where $ \psi\in C^0(\mathbb{R}^d\times \mathbb{R}_+ ,\mathbb{R}) $ is a continuous function with polynomial growth, and set $Y(t)=E[g(X_T)|{\cal F}_t]$.
%In finance, this corresponds to the case of lookback options.
Then $ g $ satisfies Assumption \ref{hyp_g1}, and our approximation method applies. Moreover, if $ \psi $ is Lipschitz-continuous, then  Theorem \ref{th_error} yields an explicit control of the approximation error  with a  $\sqrt{\log n / n}$ bound.
\end{example}

\section{Comparison with approaches based on the Malliavin calculus}\label{sec:comparison}

The vertical derivative $\nabla_X Y(t)$ which appears in the martingale representation formula may be viewed as a 'sensitivity' of the martingale $Y$ to the initial condition $X(t)$. Thus, our method is related to methods previously proposed for 'sensitivity analysis' of Wiener functionals.

One can roughly classify such methods into two categories \cite{chenglasserman07}: methods that differentiate paths and methods that differentiate densities. When the density of the functional is known, the sensitivity of an expectation with respect to some parameter is  to differentiate directly the density function with respect to the parameter. However, as this is almost never the case in a general diffusion model, let alone a non-Markovian model, alternative methods, are used: these consist of differentiating either the functional $ g $ or the process with respect to the parameter under the expectation sign, then estimating the expectation with the Monte-Carlo method. To differentiate process, one required the existence of the so-called first variation process, which requires the regularity of the coefficients of the SDE satisfied by $ X $.  

Sensitivity estimators for non-smooth functionals may be computed  using Malliavin calculus: this approach, proposed by Fourni\'e et al. \cite{FLLLT99} and developed by  Cvitanic, Ma and Zhang \cite{CMZ03}, Fournié et al. \cite{FLLL01}, Gobet and Kohatsu-Higa \cite{GK03}, Kohatsu-Higa and Montero \cite{KM04}, Davis and Johansson \cite{DJ06} and   others, uses the Malliavin  integration-by-parts formula on Wiener space in the case where $ g $ is not smooth.   These methods  require  quite demanding regularity assumptions (differentiability and ellipticity condition on $ \sigma $ for example) on the coefficients of the initial SDE satisfied by $ X $.    
    
By contrast, the approximation method  presented here allows for any continuous functional $g$ with polynomial growth and  requires only a  (Lipschitz-)continuity assumption on  the functional $ \sigma$ and allows for degenrate coefficients. It is thus applicable to a wider range of examples than the Malliavin approach, while being arguably simpler from a computational viewpoint. 
Contrarily to the Malliavin approach which involves differentiating in the Malliavin sense, then discretizing the tangent process, 
our method involves  discretizing then differentiating (the Euler scheme) which, as  argued in  \cite{chenglasserman07}, has its computational advantages.

In our setting, we have $ F_n \in \mathbb{C}^{1, 2}_{loc, r} (\Lambda_T) $ which is sufficient for   obtaining an approximation of martingale representations via the functional Itô formula. One can ask if the Euler approximation  $ {_n X} $ can also be used to obtain a Clark-Haussmann-Ocone type formula, and in this case, whether the pathwise vertical derivative $ \nabla_{\omega} F_n (t, W_t) $ leads to the same representation as the Clark-Haussmann-Ocone formula.

For $ n \in \mathbb{N},$ define $ H_n = g ( {_n X}_T (W_T) ) $ with $ {_n X} $ the weak piecewise constant Euler-Maruyama scheme defined by \eqref{euler}. By the definition of $ {_n X} $, the random variable $ H_n $ actually depends only on a finite number of Gaussian variables: $ W (t_1) $, $ W (t_2) - W (t_1) $, $ \cdots $, $ W (t_n) - W (t_{n-1}) $, thus it can be written as:
\[ H_n = h_n (W (t_1), W (t_2) - W (t_1), \cdots, W (t_n) - W (t_{n-1})) \]
with $ h_n: \mathcal{M}_{d, n} (\mathbb{R}) \rightarrow \mathbb{R} $ ($ h_n $ is actually $ g \circ p $ with $ p $ defined by \eqref{p}).

Clearly if $ h_n $ is a smooth function with polynomial growth, then $ H_n \in \mathbf{D}^{1, 2}$ with Malliavin derivative \cite{N09}:
\[ \mathbb{D}_t H_n = ( \mathbb{D}_t^k H_n, 1 \leq k \leq d ) \in \mathbb{R}^d, \qquad {\rm where} \]
\[ \mathbb{D}_t^k H_n = \sum_{j=0}^{n-1} \partial_{kj} h_n (W (t_1), W (t_2) - W (t_1), \cdots, W (t_n) - W (t_{n-1})) {\mathbf 1}_{[t_j, t_{j+1})} (t), \quad t \in [0, T) \]
In this case, assume that $ t \in [t_j, t_{j+1}) $ for some $ 0 \leq j \leq n-1 $, we have: for $ 1 \leq k \leq d $,
\begin{eqnarray*}
&& \mathbb{E} \left[ \mathbb{D}_t^k H_n | \mathcal{F}_t \right] \\
&=& \mathbb{E} \left[ \partial_{kj} h_n (W (t_1), W (t_2) - W (t_1), \cdots, W (t_n) - W (t_{n-1})) | \mathcal{F}_t \right] \\
&=& \mathbb{E} \big[ \partial_{kj} h_n (\omega (t_1), \omega (t_2) - \omega (t_1), \cdots, \omega (t_j) - \omega (t_{j-1}), W (t_{j+1}) - W (t) + \omega (t) - \omega (t_j),  \\
&& W (t_{j+2}) - W (t_{j+1}), \cdots, W (t_n) - W (t_{n-1})) \big] |_{\omega_t = W_t} \\
&=& \mathbb{E} \Big[ \frac{\partial}{\partial h} h_n (\omega (t_1), \omega (t_2) - \omega (t_1), \cdots, \omega (t_j) - \omega (t_{j-1}), W (t_{j+1}) - W (t) + \omega (t) - \omega (t_j) + h e_k, \\
&& W (t_{j+2}) - W (t_{j+1}), \cdots, W (t_n) - W (t_{n-1})) |_{h = 0} \Big] |_{\omega_t = W_t} \\
&=& \frac{\partial}{\partial h} \Big( \mathbb{E} \big[  h_n (\omega (t_1), \omega (t_2) - \omega (t_1), \cdots, \omega (t_j) - \omega (t_{j-1}), W (t_{j+1}) - W (t) + \omega (t) - \omega (t_j) + h e_k, \\
&& W (t_{j+2}) - W (t_{j+1}), \cdots, W (t_n) - W (t_{n-1})) \big]_{\omega_t = W_t} \Big) |_{h = 0}
\end{eqnarray*}
which is none other  than
\[ \partial_k F_n (t, W_t) = \lim\limits_{h \to 0} \frac{F_n (t, W_t + h e_k {\mathbf 1}_{[t, T]}) - F_n (t, W_t)}{h}. \]
So in the case where $ h_n $ are smooth, our method provides the same result as given by the Clark-Haussmann-Ocone formula applied to $h_n$.  
 However, in our framework, as the functional $ g $ is only assumed to be continuous with polynomial growth, the function $ h_n $ may fail to be differentiable. So, even in the cylindrical case, it is not clear whether the random variable $ H_n $ is differentiable in the Malliavin sense, and even if it is the case, it is difficult to obtain an explicit form for $ \mathbb{E} \left[ \mathbb{D}_t H_n | \mathcal{F}_t \right] $ using the Malliavin calculus.

However, our approximation method applies even in the cases where $ H_n $ is not differentiable in the Malliavin sense: indeed, as shown in Section 3, as soon as $H$ is square-integrable, the martingale $ {_n Y} (t) = \mathbb{E} \left[ H_n | \mathcal{F}_t \right] $ has a smooth functional representation which  is differentiable in the pathwise sense, even though $ H_n $ is not differentiable, neither in a pathwise nor in the Malliavin sense.

%% file: FIC.tex
The Functional It\^o calculus \cite{Cont2012} is a functional calculus which extends the It\^o calculus to path-dependent functionals of stochastic processes. It was first introduced in a pathwise setting  \cite{CF10A,CF10B,D09} using a notion of  pathwise derivative for functionals on the space of right-continuous functions with left limits,
and extended in \cite{CF13} to a weak calculus applicable to all square-integrable martingales, which has a natural connection to the martingale representation theorem. We recall here some key concepts and results of this approach, following \cite{Cont2012}.

%This extension defines a weak derivative which is shown to be the inverse of the It\^o integral, and provides thus a natural candidate for the integrand in the martingale representation formula. However, as a weak derivative, this quantity is not explicitly computable in general and its computation requires numerical approximations, which are the object of this paper.

 Let $ X $ be the canonical process on $ \Omega = D ( [0,T], \mathbb{R}^d ) $, and $ ( \mathcal{F}_t^0 )_{t \in [0, T]} $ be the filtration generated by $ X $. We are interested in {\it non-anticipative functionals} of $X$, that is, functionals  $ F:[0, T] \times D( [0, T], \mathbb{R}^d ) \mapsto \mathbb{R} $ such that
\begin{equation}  \forall \omega \in \Omega,\qquad F (t, \omega) = F (t, \omega_t). \label{eq:nonanticipative}\end{equation}
The process $ t \mapsto F (t,X_t) $  then only depends on the path of $X$ up to $ t $ and is $ ( \mathcal{F}_t^0 ) $-adapted.

It is convenient to define such functionals  on the space of stopped paths \cite{Cont2012}: a stopped path is an equivalence class in $ [0, T] \times D( [0, T], \mathbb{R}^d ) $ for the following equivalence relation:
\begin{equation} \label{equiv}
(t, \omega) \sim (t', \omega') \Longleftrightarrow ( t=t' \quad \mathrm{and} \quad \omega_t = \omega'_{t'} ).
\end{equation}
The space of stopped paths is defined as the quotient of $ [0, T] \times D( [0, T], \mathbb{R}^d ) $ by the equivalence relation \eqref{equiv}:
\[ \Lambda_T = \{ (t, \omega (t \wedge \cdot)), (t, \omega) \in [0, T] \times D( [0, T],\mathbb{R}^d ) \} = \left( [0, T] \times D( [0, T], \mathbb{R}^d ) \right) / \sim \]
We denote ${\cal W}_T$ the subset of $\Lambda_T$ consisting of continuous stopped paths.
We endow this set with a metric space structure by defining the following distance:
\begin{eqnarray*}
d_{\infty} ((t, \omega), (t', \omega')) = \sup_{u \in [0, T]} | \omega (u \wedge t) - \omega' (u \wedge t') | + | t - t' |
= \| \omega_t - \omega'_{t'} \|_{\infty} + | t - t' |
\end{eqnarray*}
$ ( \Lambda_T, d_{\infty} ) $ is then a complete metric space.
Any functional verifying the non-anticipativeness condition \eqref{eq:nonanticipative} can be equivalently viewed as a functional on $F:\Lambda_T\to \mathbb{R}$:
\begin{defi} \label{def_nonanti}
A non-anticipative functional on $ D( [0, T], \mathbb{R}^d ) $ is a measurable map\\ $ F $: $ (\Lambda_T, d_{\infty}) \longrightarrow \mathbb{R} $ on the space $ (\Lambda_T, d_{\infty}) $ of stopped paths.
\end{defi}
Using the metric structure of $(\Lambda_T, d_{\infty})$, we denote by $ \mathbb{C}^{0, 0} (\Lambda_T) $ the set of  continuous maps $F:(\Lambda_T, d_{\infty})\mapsto\mathbb{R}$. Some weaker notions of continuity for non-anticipative functionals turn out to be useful \cite{CF10A}:
\begin{defi} \label{def_functional}
A non-anticipative functional $ F $ is said to be:
\begin{itemize}
\item continuous at fixed times if for any $ t \in [0, T] $, $ F (t, \cdot) $ is continuous with respect to the uniform norm $ \| \cdot \|_{\infty} $ in $ [0, T] $, i.e. $ \forall \omega \in D( [0, T], \mathbb{R}^d ) $, $ \forall \epsilon > 0 $, $ \exists \eta > 0 $, $ \forall \omega' \in D( [0, T], \mathbb{R}^d ) $,
\[  \| \omega_t - \omega'_t \| < \eta \implies | F (t, \omega) - F (t, \omega') | < \epsilon \]
\item left-continuous if $ \forall (t, \omega) \in \Lambda_T $, $ \forall \epsilon > 0 $, $ \exists \eta > 0 $ such that $ \forall (t', \omega') \in \Lambda_T $,
\[ ( t' < t \quad \mathrm{and} \quad d_{\infty} \left( (t, \omega), (t', \omega') \right) < \eta ) \implies | F (t, \omega) - F (t', \omega') | < \epsilon \]
We denote by $ \mathbb{C}^{0, 0}_l (\Lambda_T) $ the set of left-continuous functionals. Similarly, we can define the set $ \mathbb{C}^{0, 0}_r (\Lambda_T) $ of right-continuous functionals.
\end{itemize}
\end{defi}
We also introduce a notion of local boundedness for functionals.
\begin{defi} \label{def_bound}
A non-anticipative functional $ F $ is said to be boundedness-preserving if for every compact subset $ K $ of $ \mathbb{R}^d $, $ \forall t_0 \in [0, T] $, $ \exists C (K, t_0) > 0 $ such that:
\[ \forall t \in [0, t_0],\quad \forall (t, \omega) \in \Lambda_T, \quad\omega ([0, t]) \subset K \implies F (t, \omega) < C (K, t_0). \]
We denote by $ \mathbb{B} (\Lambda_T) $ the set of boundedness-preserving functionals.
\end{defi}
We now recall some notions of differentiability for functionals following \cite{CF13,Cont2012}. For $ e \in \mathbb{R}^d $ and $ \omega \in D ( [0,T], \mathbb{R}^d ) $, we define the vertical perturbation $ \omega_t^e $ of $(t,\omega)$ as the c\`adl\`ag path obtained by shifting the path by $e$ after $ t $:
$$\omega_t^e = \omega_t + e {\mathbf 1}_{[t, T]}.$$
\begin{defi} \label{def_diff}
A non-anticipative functional $ F $ is said to be:
\begin{itemize}
\item horizontally differentiable at $ (t, \omega) \in \Lambda_T $ if
\[ \mathcal{D} F (t, \omega) = \lim\limits_{h \to 0^+} \frac{F (t+h, \omega) - F (t, \omega)}{h} \]
exists. If $\mathcal{D} F (t, \omega)$ exists for all $ (t, \omega) \in \Lambda_T $, then the non-anticipative functional $ \mathcal{D} F $ is called the horizontal derivative of $ F $.
\item vertically differentiable at $ (t, \omega) \in \Lambda_T $ if the map:
\[
\begin{array}{ccc}
\mathbb{R}^d & \longrightarrow & \mathbb{R} \\
e & \mapsto & F (t, \omega_t + e {\mathbf 1}_{[t, T]})
\end{array}
 \]
is differentiable at  $ 0 $. Its gradient at $ 0 $ is called the vertical derivative of $ F $ at $ (t, \omega) $:
\[ \nabla_{\omega} F (t, \omega) = ( \partial_i F (t, \omega), i = 1, \cdots, d) \in \mathbb{R}^d \]
with
\[ \partial_i F (t, \omega) = \lim_{h \to 0} \frac{F (t, \omega_t + h e_i {\mathbf 1}_{[t, T]}) - F (t, \omega_t)}{h} \]
where $ ( e_i, i=1, \cdots, d ) $ is the canonical basis of $ \mathbb{R}^d $. If $ F $ is vertically differentiable at all $ (t, \omega) \in \Lambda_T $,   $ \nabla_{\omega} F: (t,\omega)\to \mathbb{R}^d$ defines a non-anticipative map called the vertical derivative of $ F $.
\end{itemize}
\end{defi}
We may repeat the same operation on $ \nabla_{\omega} F $ and define similarly $ \nabla_{\omega}^2 F $, $ \nabla_{\omega}^3 F $, $ \cdots $.
This leads us to define the the following classes of smooth functionals:
\begin{defi}[Smooth functionals] \label{def:C12}
We define
$ \mathbb{C}^{1, k}_b (\Lambda_T) $ as the set of non-anticipative functionals
 $ F: (\Lambda_T,d) \to \mathbb{R}$ which are
\begin{itemize}
\item horizontally differentiable with $ \mathcal{D} F $ continuous at fixed times;
\item $ k $ times vertically differentiable with $ \nabla_{\omega}^j F \in \mathbb{C}^{0, 0}_l (\Lambda_T) $ for $ j = 0, \cdots, k $;
\item  $ \mathcal{D} F, \nabla_{\omega} F, \cdots, \nabla_{\omega}^k F \in \mathbb{B} (\Lambda_T) $.
\end{itemize}
\end{defi}
We denote $ \mathbb{C}^{1, \infty}_b(\Lambda_T) = \cap_{k \geq 1} \mathbb{C}^{1, k}_b (\Lambda_T).$

Many examples of functionals may fail to be globally smooth, but their derivatives may still be well behaved except at certain stopping times, which motivates the following definition \cite{Cont2012}:
\begin{defi} \label{def_loc}
A non-anticipative functional $ F $ is said to be locally regular of class $ \mathbb{C}_{\rm loc}^{1, 2} (\Lambda_T) $ if there exists an increasing sequence $ (\tau_n)_{n \geq 0} $ of stopping times with $ \tau_0 = 0 $ and $ \tau_n \underset{n \to \infty}{\longrightarrow} \infty $, and a sequence of functionals $ F_n \in \mathbb{C}_b^{1, 2} (\Lambda_T) $ such that:
\[ F (t, \omega) = \sum_{n \geq 0} F_n (t, \omega) {\mathbf 1}_{[\tau_n (\omega), \tau_{n+1} (\omega))} (t), \quad \forall (t, \omega) \in \Lambda_T \]
\end{defi}
We recall now the functional It\^o formula for non-anticipative functionals of a continuous semimartingale \cite[Theorem 4.1]{CF13}:
\begin{prop}[\cite{CF10A,CF13}]  \label{prop_ito}
Let $S$ be a continuous semimartingale defined on a probability space $ ( \Omega, \mathcal{F}, \mathbb{P} ) $. For any non-anticipative functional $ F \in \mathbb{C}^{1, 2}_{\rm loc} (\Lambda_T) $ and any $ t \in [0, T] $, we have:
\begin{eqnarray*}
F (t, S_t) - F (0, S_0) &=& \int_0^t \mathcal{D} F (u, S_u) du + \int_0^t \nabla_{\omega} F (u, S_u) \cdot dS(u)
+ \frac{1}{2} \int_0^t \mathrm{tr} \left( \nabla_{\omega}^2 F (u, S_u) d [S] (u) \right)
%\quad \mathbb{P} \mathrm{-a.s.}
\end{eqnarray*}
\end{prop}

Actually the same functional It\^o formula may also be obtained for functionals whose vertical derivatives are right-continuous rather than left-continuous. We denote by $ \mathbb{C}^{1, 2}_{b, r} (\Lambda_T) $ the set of non-anticipative functionals $ F $ satisfying:
\begin{itemize}
\item
$ F $ is horizontally differentiable with $ \mathcal{D} F $ continuous at fixed times;
\item
$ F $ is twice vertically differentiable with $ F \in \mathbb{C}_l^{0, 0} (\Lambda_T) $ and $ \nabla_{\omega} F, \nabla^2_{\omega} F \in \mathbb{C}_r^{0, 0} (\Lambda_T) $;
\item
$ \mathcal{D} F, \nabla_{\omega} F, \nabla^2_{\omega} F \in \mathbb{B} (\Lambda_T) $;
\end{itemize}
The localization is more delicate in this case, and we are not able to state a local version of the functional It\^o formula by simply replacing $ F_n \in \mathbb{C}_b^{1, 2} (\Lambda_T) $ by $ F_n \in \mathbb{C}_{b, r}^{1, 2} (\Lambda_T) $ in \textbf{Definition \ref{def_loc}} (see Remark $ 4.2 $ in \cite{F10}). However if the stopping times $ \tau_n $ are  deterministic, then the functional It\^o formula is still valid (Proposition $ 2.4 $ and Remark $ 4.2 $ in \cite{F10}).
\begin{defi} \label{def_loc2}
A non-anticipative functional is said to be locally regular of class $ \mathbb{C}_{\rm loc, r}^{1, 2} (\Lambda_T) $ if there exists an increasing sequence $ (t_n)_{n \geq 0} $ of deterministic times with $ t_0 = 0 $ and $ t_n \underset{n \to \infty}{\longrightarrow} \infty $, and a sequence of functionals $ F_n \in \mathbb{C}_{b, r}^{1, 2} (\Lambda_T) $ such that:
\[ F (t, \omega) = \sum_{n \geq 0} F_n (t, \omega) {\mathbf 1}_{[t_n, t_{n+1})} (t), \quad \forall (t, \omega) \in \Lambda_T \]
\end{defi}

\begin{prop}[\cite{CF13}] \label{prop_ito2}
Let $S$ be a continuous semimartingale defined on a probability space $ ( \Omega, \mathcal{F}, \mathbb{P} ) $. For any non-anticipative functional $ F \in \mathbb{C}^{1, 2}_{loc, r} (\Lambda_T) $ and any $ t \in [0, T] $, we have:
\begin{eqnarray*}
F (t, S_t) - F (0, S_0) &=& \int_0^t \mathcal{D} F (u, S_u) du + \int_0^t \nabla_{\omega} F (u, S_u) \cdot dS(u) \\
&&+ \frac{1}{2} \int_0^t \mathrm{tr} \left( \nabla_{\omega}^2 F (u, S_u) d [S] (u) \right) \quad \mathbb{P} \mathrm{-a.s.}
\end{eqnarray*}
\end{prop}
Finally we present briefly the martingale representation formula established in \cite{CF13}. Let $ (X_t)_{t \in [0, T]} $ be a continuous $ \mathbb{R}^d $-valued martingale defined on a probability space $ ( \Omega, \mathcal{F}, \mathbb{P} ) $ with absolutely continuous quadratic variation:
\[ [ X ] (t) = \int_0^t A (u) du \]
where $ A $ is a $ \mathcal{M}_d (\mathbb{R}) $-valued process.
Denote by $ (\mathcal{F}_t^X) $ the natural filtration of $ X $ and $ \mathcal{C}_b^{1, 2} (X) $ the set of $ ( \mathcal{F}_t^X ) $-adapted processes $ Y $ which admit a functional representation in $ \mathbb{C}_b^{1, 2} (\Lambda_T) $:
\begin{equation} \label{repre}
\mathcal{C}_b^{1, 2} (X) = \{ Y, \exists F \in \mathbb{C}_b^{1, 2} (\Lambda_T), Y(t) = F (t, X_t) \quad dt \times d \mathbb{P} \mathrm{-a.e.} \}
\end{equation}
If $ A(t) $ is non-singular almost everywhere, i.e. $ \mathrm{det} (A(t)) \neq 0 $, $ dt \times d \mathbb{P} $-a.e., then for any $ Y \in \mathcal{C}_b^{1, 2} (X) $, the predictable process
\[ \nabla_X Y (t) = \nabla_{\omega} F (t, X_t) \]
is uniquely defined up to an evanescent set, independently of the choice of $ F \in \mathbb{C}_b^{1, 2} (\Lambda_T) $ in the representation \eqref{repre}.
This process $ \nabla_X Y $ is called the vertical derivative of $ Y $ with respect to $ X$.
For martingales which are smooth functionals of $X$, the operator  $\nabla_X: \mathcal{C}_b^{1, 2} (X) \mapsto \mathcal{C}^{0, 0}_l(X)$ yields the integrand in the martingale representation theorem:
\begin{cor} \label{cor_repre}
If $ Y \in \mathcal{C}_b^{1, 2} (X) $ is a square-integrable martingale, then
\[ \forall t \in [0, T], \quad Y (t) = Y (0) + \int_0^t \nabla_X Y \cdot dX \quad \mathbb{P} \mathrm{-a.s.} \]
\end{cor}
Consider now the case where $X$ is a square-integrable martingale. Let $ \mathcal{M}^2 (X) $ be the space of square-integrable $ ( \mathcal{F}_t^X) $-martingales with initial value zero, equipped with the norm $ \| Y \|^2 = \sqrt{\mathbb{E} | Y(T) |^2} $. Cont \& Fourni\'e \cite[Theorem 5.8]{CF13} show that the operator  $\nabla_X: \mathcal{C}_b^{1, 2} (X) \mapsto \mathcal{C}^{0, 0}_l(X)$ admits a unique continuous extension to a weak derivative  $ \nabla_X:\mathcal{M}^2 (X) \to \mathcal{L}^2 (X)$ which satisfies
 the following martingale representation formula:
\begin{prop}[\cite{CF13}] \label{prop_mart}
For any square-integrable $ ( \mathcal{F}_t^W ) $-martingale $ Y $, we have:
\[ \forall t \in [0, T], \quad Y(t) = Y(0) + \int_0^t \nabla_X Y \cdot dX \quad \mathbb{P} \mathrm{-a.s.} \]
\end{prop}
This weak vertical derivative $ \nabla_X Y $ coincides with the pathwise vertical derivative $ \nabla_{\omega} F (t, X_t) $  when $ Y $ admits a locally regular functional representation,  i.e. $ Y (t) = F (t, X_t) $ with $ F \in \mathbb{C}_{loc}^{1, 2} (\Lambda_T) \cup \mathbb{C}_{loc, r}^{1, 2} (\Lambda_T) $. For a general square-integrable martingale $ Y $, the weak derivative $ \nabla_X Y $ is not directly computable through a pathwise perturbation. An approximation procedure is thus necessary for computing $\nabla_XY$. The result of  \cite{CF13} guarantees the existence of such approximations; in the sequel we propose explicit, and computable, constructions of such approximations.

%% file: SPA2851.bbl
\def\cprime{$'$}